\documentclass[10pt, a4paper]{amsart}
\usepackage{authblk}
\usepackage{bbm}
\usepackage{amssymb}
\usepackage{mathtext}
\usepackage{amsmath}
\usepackage{amsfonts}
\usepackage{amssymb}
\usepackage{amsthm}
\usepackage{mathrsfs}
\usepackage{hyperref}
\usepackage{amsaddr}
\usepackage{graphicx}
\usepackage{geometry} 
\newtheorem{theorem}{Theorem}[section]
\newtheorem{assumption}{Assumption}
\newtheorem{proposition}{Proposition}[section]
\newtheorem{corollary}{Corollary}[section]
\newtheorem{lemma}{Lemma}[section]
\newtheorem{definition}{Definition}[section]

\newtheorem{remark}{Remark}[section]

\newcommand{\E}{\ensuremath{\mathbb E}}
\newcommand*{\QED}{\hfill\ensuremath{\square}}

\begin{document}

\markboth{Yu. Mishura, A. Yurchenko-Tytarenko}
{Option pricing in fractional Heston-type model}


\title{
OPTION PRICING IN FRACTIONAL HESTON-TYPE MODEL
}

\author{YULIYA MISHURA}
\address{Faculty of Mechanics and Mathematics, Taras Shevchenko National University of Kyiv, Akad. Glushkova Av. 4-e, Kyiv, 03127, Ukraine}

\author{ANTON YURCHENKO-TYTARENKO}
\address{Faculty of Mechanics and Mathematics, Taras Shevchenko National University of Kyiv, Akad. Glushkova Av. 4-e, Kyiv, 03127, Ukraine}
\email{antonyurty@gmail.com}

\maketitle

\begin{abstract}
In this paper, we consider option pricing in a framework of the fractional Heston-type model with $H>1/2$. As it is impossible to obtain an explicit formula for the expectation $\E f(S_T)$ in this case, where $S_T$ is the asset price at maturity time and $f$ is a payoff function, we provide a discretization schemes $\hat Y^n$ and $\hat S^n$ for volatility and price processes correspondingly and study convergence  $\E f(\hat S^n_T) \to \E f(S_T)$ as the mesh of the partition tends to zero. The rate of convergence is calculated. As we allow $f$ to have discontinuities of the first kind which can cause errors in straightforward Monte-Carlo estimation of the expectation, we use Malliavin calculus techniques to provide an alternative formula for $\E f(S_T)$ with smooth functional under the expectation.
\end{abstract}


\section{Introduction}

Despite its undoubtedly significant historical and theoretical value, the classical Black-Scholes model does not explain numerous empirical phenomena that can be observed on real-life markets, such as implied volatility smile and skew. In order to overcome this issue, \cite{HW1987} and, later, \cite{Hest} introduced stochastic volatility models that emerged into an essential subject of research activity in financial modeling nowadays. 

To illustrate the range of existing models (without trying to list all possible references), we recall the approaches of \cite{AltNeu}, \cite{BNS1}, \cite{BNS2}, \cite{CC2002}, \cite{CT2004}, \cite{KuprSchou2005}, \cite{NicVen}, \cite{Shep1996},  and so on.

A separate class of stochastic volatility models are those based on fractional Brownian motion. They allow to reflect the so-called ``memory phenomenon'' of the market (for more detail on market models with memory see, for instance, \cite{AI2005, Ding1993, Yam2005}). In this context, we should also mention \cite{BorMik1996, ChronViens2012, Comte2012} and \cite{BdPM}.   

In the present paper, we consider option pricing in a framework of the fractional modification of the Heston-type model, namely a financial market with a finite maturity time $T$ that is composed of two assets:

(i) a risk-free bond (or bank account) $B=\{B_t,~t\in[0,T]\}$, the dynamics of which is characterized by the formula
\begin{equation}\label{fHestonRiskFree}
B_t = e^{\lambda t}, \quad t\in [0,T],
\end{equation}
where $\lambda \in \mathbb R^+$ represents the risk-free interest rate;
 
(ii) a risky asset $S = \{S_t,~t\in[0,T]\}$, the evolution in time of which is given by the system of stochastic differential equations
\begin{equation}\label{fHestonPrice}
dS_t = \mu S_tdt + \sigma(Y_t)S_tdW_t,
\end{equation}
\begin{equation}\label{fHestonVolatility}
dY_t = \frac{1}{2}\left(\frac{\kappa}{Y_t} -\theta Y_t\right)dt +\frac{\nu}{2} dB_t^H, \quad t\in[0,T],
\end{equation}
with non-random initial values $S_0, Y_0>0$, where the process $W=\{W_t,~t\ge 0\}$ is a standard Wiener process, $\mu \in\mathbb R$, $\kappa, \theta, \nu >0$ are constants, $\sigma$: $[0,\infty) \to [0, \infty)$ is a function that satisfies some regularity properties and $B^H=\{B^H_t, ~t\in[0,T]\}$ is a fractional Brownian motion with the Hurst index $\frac{1}{2}<H<1$, which corresponds to the ``long memory'' case. $W$ and $B^H$ are assumed to be correlated.

The process $Y$ was extensively studied in \cite{MYuT1, MYuT2} and, for the case $\kappa=0$, in \cite{MPRYuT}. Note that, according to \cite{NualOuk}, the process $Y$ exists, is unique and has continuous paths until the first moment of zero hitting. Moreover, in Theorem 2 of \cite{MYuT1} it was shown that in case of $\kappa>0$ and $H>\frac{1}{2}$ such process is strictly positive and never hits zero, therefore exists, is unique and continuous on the entire $[0,T]$. 

Such choice of the volatility process can be explained by the fact that $Y$ can be interpreted as the square root of the fractional version of Cox-Ingersoll-Ross process. Indeed, according to \cite{MYuT1}, Theorem 1, the process $X=\{Y^2(t),~t\in[0,T]\}$ satisfies the stochastic differential equation of the form 
\begin{equation*}
dX_t = (\kappa-\theta X_t)dt +\nu \sqrt{X_t} dB_t^H, \quad X_0=Y_0^2>0,
\end{equation*}
until the first moment of zero hitting, where the integral $\int_0^t \sqrt{X_s} dB_s^H$ is considered as the pathwise limit of the sums
\begin{equation*}
\sum_{k=1}^n \frac{X_{t_k}+X_{t_{k-1}}}{2}(B^H_{t_k}-B^H_{t_{k-1}}),
\end{equation*}
as the mesh of the partition $0 = t_0 < t_1 < ... < t_n =t$ tends to zero.

Note that, due to Kolmogorov theorem, fractional Brownian motion $B^H$ has a modification with H\"older continuous paths up to order $H$. Hence, from the form of the equation \eqref{fHestonVolatility}, the process $Y$ also has a modification with trajectories that are H\"older-continuous up to order $H$. Therefore, in case of $H>\frac{1}{2}$, the sum of H\"older exponents of the integrator and integrand in the integral
\begin{equation*}
\int_0^t \sqrt{X_s} dB_s^H = \int_0^t Y_s  dB_s^H
\end{equation*}
exceeds 1 and, due to \cite{Zahle}, the corresponding integral exists as the pathwise limit of Riemann-Stieltjes integral sums.

It should be also mentioned that for the case $H<1/2$, the process $Y$ can hit zero and it is not clear whether the solution exists on the entire $[0,T]$ (see \cite{MYuT2} for more detail). Therefore, we will concentrate on the case $H>1/2$. For more information on markets with rough volatility see, for example, \cite{GatJaiRos} or \cite{NeuRos}.

An analogue of the model \eqref{fHestonPrice}, \eqref{fHestonVolatility} was considered in \cite{BdPM} with fractional Ornstein-Uhlenbeck process instead of $Y$. However, Ornstein-Uhlenbeck process can take negative values with positive probability which is a notable drawback for a stochastic volatility model. 

Note that it is impossible to calculate $\mathbb E f(S_T)$ (with $f$ being a payoff function) for option pricing analytically, so numerical methods should be used. Therefore it is required to provide a decent discretization scheme for $S_T$ and prove the convergence 
\begin{equation}\label{intro: convergence}
\mathbb E f(\hat S^n_T) \to \mathbb E f(S_T), \qquad n\to\infty,
\end{equation}
where $\hat S^n$ is a discretized version of the process $S$. Moreover, we allow $f$ to have discontinuities of the first kind which can cause errors in straightforward Monte-Carlo estimation of the expectation, so we provide an alternative formula with smooth functional under the expectation. In such framework, we also give the rate of convergence \eqref{intro: convergence}.

It should be mentioned that the market with risky asset defined by \eqref{fHestonPrice}--\eqref{fHestonVolatility} is arbitrage-free, incomplete but admits minimal martingale measure (see Section \ref{sec: model properties}). However, the expectations calculated with respect to the minimal martingale and objective measures differ only by non-random coefficient, therefore, for simplicity, we concentrate on expectation with respect to the objective measure. In order to model the volatility $Y$, we use the inverse Euler approximation scheme studied in \cite{HHKW}.

The paper is organized as follows.  In Section \ref{sec: model description}, we describe main assumptions concerning relation between the Wiener process and the fractional Brownian motion as well as volatility function $\sigma$ and payoff function $f$. In Section \ref{sec: model properties} several important properties of both price and volatility processes are presented and the arbitrage-free property is discussed. In Section \ref{sec: option pricing} we apply the Malliavin calculus techniques, following \cite{AltNeu} and \cite{BdPM}, to obtain the formula for option price that does not contain discontinuities (which are allowed for the payoff function $f$). In Section \ref{sec: Inverse Euler}, we study the rate of convergence of Monte-Carlo estimation of the option price $\E f(S_T)$ based on inverse Euler approximation scheme for fractional CIR process presented in \cite{HHKW}. In Section \ref{sec: simulations}, we give results of numerical simulations for different payoff functions $f$. Section \ref{sec: proofs} contains the proofs of all results of the paper. Appendix \ref{sec: preliminaries} is devoted to several well-known results from the Malliavin calculus used in this paper.

\section{Model description and main assumptions}\label{sec: model description}

Consider the market with risk-free asset $B$ given by \eqref{fHestonRiskFree} and risky asset $S$, the dynamics of which is described by stochastic differential equations \eqref{fHestonPrice}, \eqref{fHestonVolatility}.

Denote
\begin{equation*}
K(t,s) = c_Hs^{\frac{1}{2}-H}\int_s^t u^{H-\frac{1}{2}}(u-s)^{H-\frac{3}{2}}du\mathbbm{1}_{s<t},
\end{equation*}
\begin{equation*}
c_H = \left(\frac{H(2H-1)}{B(2-2H, H-\frac{1}{2})}\right)^{1/2},
\end{equation*}
where $B(\cdot,\cdot)$ is the Beta function. Then, according to \cite{Norros1999}, the process $B^H = \{B^H_t,~t\in[0,T]\}$ given by 
\begin{equation}\label{kernelRepresentation}
B^H_t = \int_0^t K(t,s)dV_s, \quad t\in[0,T],
\end{equation}
where $V = \{V_t,~t\in[0,T]\}$ is a Wiener process, is the fractional Brownian motion with Hurst parameter $H$.

The processes $W$ and $B^H$ from \eqref{fHestonPrice}, \eqref{fHestonVolatility} are assumed to be correlated and the form of the dependence is defined on the basis of representation \eqref{kernelRepresentation} as follows. 

\begin{assumption}\label{AssumCorr}
The processes $W$ and $V$ from \eqref{fHestonPrice} and \eqref{kernelRepresentation} correspondingly are correlated:
\begin{equation*}
\mathbb E W_t V_t = \rho t, \quad t\in [0,T],
\end{equation*}
with some constant $\rho\in [-1,1]$.
\end{assumption}

\begin{remark}\label{two Wiener processes}
Assumption \ref{AssumCorr} means that $W_t = \rho V_t + \sqrt{1-\rho^2}\tilde {V}_t$, $t\in[0,T]$, where $\tilde V$ is a Wiener process independent of $V$.
\end{remark}

The function $\sigma$: $\mathbb{R}\to\mathbb{R}$ is assumed to satisfy the following conditions.

\begin{assumption}\label{AssumSigma}
For some constant $C_\sigma>0$:
\begin{itemize}
\item[(i)] there exists such $\sigma_{\min}>0$ that for all $x\in\mathbb{R}$: $\sigma(x)>\sigma_{\min}>0$;
\item[(ii)] $\sigma$ has moderate polyniomial growth, i.e. there is such $q \in (0,1)$ that 
\begin{equation*}
\sigma(x)\le C_\sigma (1+|x|^q), \quad x\in \mathbb R;
\end{equation*}
\item[(iii)] $\sigma$ is uniformly H\"older continuous, i.e. there is such $r\in(0,1]$ that
\begin{equation*}
|\sigma(x)-\sigma(y)| \le C_\sigma |x-y|^r, \quad x,y\in\mathbb R;
\end{equation*}
\item[(iv)] $\sigma$ is differentiable a.e. w.r.t. the Lebesgue measure on $\mathbb{R}$ and there exists such $q'>0$ that
\begin{equation*}
\sigma'(x)\le C_\sigma (1+|x|^{q'}) \quad a.e.
\end{equation*}
\end{itemize}
\end{assumption}

\begin{remark}\label{remark on the model} 1) Item (i) in Assumption \ref{AssumSigma} is required for theoretical calculations as we will divide on $\sigma$ in what follows.

2) Item (ii) is necessary to ensure the finiteness of expectations of the form
\begin{equation*}
\mathbb{E}\left[\exp\left\{\int_0^t \sigma(Y_s)dW_s \right\}\right]
\end{equation*}
in case if the Wiener process $W$ and the fractional Brownian motion $B^H$ from \eqref{fHestonPrice} and \eqref{fHestonVolatility} are correlated (see Remark \ref{rem: independent case} for discussion). Note that in standard Heston model moment explosions may appear as well, see e.g. \cite{AndPit}.

3) (ii) follows from (iii) in the case $r<1$, while in (iii) we also allow $r=1$.
\end{remark}

In the framework above, we consider an option with a measurable payoff function $f:~\mathbb{R}^+\to\mathbb{R}^+$ depending on the value $S_T$ of the stock at maturity time $T$ which satisfies the following properties:

\begin{assumption}\label{AssumPayoff}
For some constant $C_\sigma>0$:
\begin{itemize}
\item[(i)] $f$ is of polynomial growth, i.e. there are such $C_f>0$ and $p>0$ that
\begin{equation*}
f(x) \le C_f (1+x^p).
\end{equation*}
\item[(ii)] $f$ is locally Riemann integrable, possibly, having discontinuities of the first kind.
\end{itemize}
\end{assumption}

\begin{remark}
In what follows, we will denote $C$ any positive constant that does not depend on time variable or diameter of the partition and the exact value of which is not important. Note that $C$ may change from line to line (and even within one line).
\end{remark}

\section{Model properties}\label{sec: model properties}

\subsection{Properties of stochastic volatility process}

In what follows we will require an auxiliary result, presented in Corollary 2.2 of \cite{MYuT2}.

\begin{theorem}\label{UpperBoundForY}
For all $H\in(0,1)$, $T>0$ and $p>0$ there are such non-random constants $C_1=C_1(T,p,Y_0, \kappa, \theta)>0$ and $C_2 = (T, p, \theta, \sigma)>0$ that for all $t\in[0,T]$:
\begin{equation*}
Y_t^p \le C_1 + C_2 \sup_{s\in [0,T]} |B^H_s|^p.
\end{equation*}

Furthermore,
\begin{equation*}
\sup_{t\in[0,T]}\E Y_t^p < \infty.
\end{equation*}
\end{theorem}


The next result is crucial for obtaining discrete approximation scheme for the process $Y$ and was presented in \cite{HHKW}.

\begin{theorem}\label{theo inv moments}
Let $p>0$ and $\kappa$, $\theta$, $\nu$ and $T$ are such that for all $t\in[0,T]$:
\begin{equation}\label{condition for inverse moments}
\kappa\exp\left\{\frac{\theta t}{2}\right\} \ge H(2H-1)(p+1)\int_0^t \frac{\nu^2}{2} \exp\left\{\frac{\theta s}{2}\right\}|t-s|^{2H-2}ds.
\end{equation}

Then there is such constant $C=C(T, Y_0, \theta)$ that
\begin{equation*}
\sup_{t\in[0,T]}\E\left[\frac{1}{Y^p_t}\right] < C.
\end{equation*}
\end{theorem}

\begin{remark}
Condition \eqref{condition for inverse moments} is satisfied if, for example,
\begin{equation*}
p+1 \le \frac{2\kappa}{\nu^2 H T^{2H-1}}.
\end{equation*}
See Remarks 3.1 and 3.2 in \cite{HHKW} for discussion.
\end{remark}

Note that condition \eqref{condition for inverse moments} involves $T$ and does not guarantee the existence of the inverse moments on whole $\mathbb R_+$. However, the following result concerning the integrated inverse moments of the volatility process $Y$ holds true.

\begin{theorem}\label{th. int. inv. m.}
Let $\beta\in\left(0,\min\{1, \frac{\kappa}{\nu^2 H T^{2H-1}}\}\right)$. Then, for all $0\le t_0 < t\le T$:
\begin{equation*}
\E \left[\int_{t_0}^t \frac{1}{Y_u^{1+\beta}}du\right] \le \frac{4}{\kappa(1-\beta)}\E (Y_t^{1-\beta} -  Y_{t_0}^{1-\beta} )  + \frac{2\theta}{\kappa} \int_{t_0}^t \E Y_u^{1-\beta} du <\infty.
\end{equation*}
\end{theorem}

\begin{theorem}\label{theor increments}
Let $\beta \in \left(0,\min\{1, \frac{\kappa}{\nu^2 H T^{2H-1}}\}\right)$. Then, there is such $C=C(\kappa,\theta,\nu, T, \beta)>0$ that for any $0\le s < t \le T$:
\begin{equation*}
\mathbb E |Y_t - Y_s|^{1+\beta} \le C|t-s|^\beta.
\end{equation*}
\end{theorem}

\begin{remark}\label{increments under China condition}
Let $p>1$ and for all $t\in[0,T]$:
\begin{equation*}
\kappa\exp\left\{\frac{\theta t}{2}\right\} \ge H(2H-1)(1+p)\int_0^t \frac{\nu^2}{2} \exp\left\{\frac{\theta s}{2}\right\}|t-s|^{2H-2}ds,
\end{equation*}
 i.e., due to Theorem \ref{theo inv moments},
\begin{equation*}
\sup_{t\in[0,T]}\E\left[\frac{1}{Y^p_t}\right] < \infty.
\end{equation*}

Proceeding just as in proof of Theorem \ref{theor increments} and taking into account that 
\begin{equation*}
(t-s)^{p-1} \E\int_s^t \frac{1}{Y_u^{p}}du < C(t-s)^p,
\end{equation*}
we can easily obtain that
\begin{equation*}
\mathbb E |Y_t - Y_s|^{p} \le C|t-s|^{pH}.
\end{equation*}
\end{remark}

\subsection{Properties of the price process}

Now let us consider several properties of the price process $S$ defined by the stochastic differential equation \eqref{fHestonPrice}.
\begin{theorem}\label{solution to price equation}
\begin{itemize}
\item[1.] For any $x>0$ and $\varrho\in[0,2)$:
\begin{equation}\label{YExponentialMoments}
\mathbb E \exp\left\{x\sup_{t\in [0,T]} |Y_t|^\varrho\right\}<\infty.
\end{equation}
\item[2.] Equation \eqref{fHestonPrice} has a unique solution of the form
\begin{equation}\label{ExplicitRepresentationOfS}
S_t = S_0 \exp\left\{\mu t +\int_0^t \sigma(Y_s)dW_s - \frac{1}{2}\int_0^t \sigma^2(Y_s)ds \right\}.
\end{equation}
\end{itemize}
\end{theorem}

\begin{remark}\label{rem: independent case}
As it was mentioned in Remark \ref{remark on the model}, presence of function $\sigma$ in \eqref{fHestonPrice}, the choice of which is restricted by Assumption \ref{AssumSigma}, is required to ensure finiteness of the moments of the form
\begin{equation*}
\mathbb{E}\left[\exp\left\{\int_0^t \sigma(Y_s)dW_s \right\}\right].
\end{equation*}

Note that Assumption \ref{AssumSigma}, (i) and (ii), does not allow $\sigma$ to be linear function, i.e. we do not consider straigthforward modification of the Heston model of the form
\begin{equation}\label{fHestonPriceInd}
dS_t = \mu S_tdt + \sigma Y_t S_tdW_t,
\end{equation}
\begin{equation}\label{fHestonVolatilityInd}
dY_t = \frac{1}{2}\left(\frac{\kappa}{Y_t} -\theta Y_t\right)dt +\frac{\nu}{2} dB_t^H, \quad t\in[0,T],
\end{equation}
where $\mu\in\mathbb R$, $\kappa$, $\theta$, $\nu$, $\sigma>0$ are constants.

However, in case of independent $W$ and $B^H$, i.e. when $\rho = 0$ in Assumption \ref{AssumCorr}, it is easy to see (e.g. by conditioning on $Y$ and solving the conditioned equation) that equation \eqref{fHestonPriceInd} has a unique solution of the form
\[
S_t = S_0 \exp\left\{\mu t + \sigma \int_0^t Y_s dW_s - \frac{\sigma^2}{2} \int_0^t Y^2_s ds\right\}.
\]

Moreover, $\E S_t < \infty$ for all $t\in[0,T]$, because the process $\tilde S$, such that
\[
\tilde S_t = \exp\left\{\sigma \int_0^t Y_s dW_s - \frac{\sigma^2}{2} \int_0^t Y^2_s ds\right\},
\]
is a non-negative local martingale and, therefore, a supermartingale.
\end{remark}

\subsection{Arbitrage-free property and incompleteness}

For the market \eqref{fHestonRiskFree}--\eqref{fHestonVolatility}, we can obtain the following result which is similar to the one in \cite{BdPM}, Theorem 4.

\begin{theorem}\label{th. arb-free}
Let the function $\sigma$ satisfy Assumption \ref{AssumSigma}. Then the market \eqref{fHestonRiskFree}--\eqref{fHestonVolatility} has the following properties.
\begin{itemize}
\item[$(i)$] It is arbitrage-free and incomplete.
\item[$(ii)$] Any probability measure $\mathbb Q$ such that
\[
\frac{d \mathbb Q}{d\mathbb P} = \exp\left\{\int_0^T \eta_1(s) dV_s + \int_0^T \eta_2(s) d\tilde V_s - \frac{1}{2}\sum_{i=1}^2 \int_0^T \eta_i^2(s)ds\right\},
\]
where $\eta_i$, $i=1,2$, are non-anticipative, bounded and satisfy the condition
\[
\rho\eta_1(s) + \sqrt{1-\rho^2}\eta_2(s) = \frac{\lambda-\mu}{\sigma(Y_s)},
\]
is a martingale measure.
\item[$(iii)$] Taking $\eta_1 = \rho\frac{\lambda-\mu}{\sigma(Y_s)}$ and $\eta_2 = \sqrt{1-\rho^2}\frac{\lambda-\mu}{\sigma(Y_s)}$, we get the minimal martingale measure.
\end{itemize}
\end{theorem}


\section{Option pricing in fractional Heston model}\label{sec: option pricing}

In this section, we will use the tools of Malliavin calculus to obtain the formula that can be used for computation of
\begin{equation*}
\E f(S_T).
\end{equation*}

Consider two-dimensional Wiener process $(V,\tilde V)$, where $V$ is given in Volterra representation \eqref{kernelRepresentation} and $\tilde V$ is defined in Remark \ref{two Wiener processes}. Denote $(D^V, D^{\tilde V})$ the stochastic derivative with respect to the two-dimensional Wiener process $(V,\tilde V)$ and recall $K$ is the kernel from representation \eqref{kernelRepresentation}. Denote also
\begin{equation}\label{logS process}
\begin{aligned}
X_t := \log S_t &= \log S_0 + \mu t - \frac{1}{2}\int_0^t \sigma^2(Y_s)ds + \int_0^t \sigma(Y_s)dW_s 
\\
&= \log S_0 + \mu t - \frac{1}{2}\int_0^t \sigma^2(Y_s)ds + \rho\int_0^t \sigma(Y_s)dV_s + \sqrt{1-\rho^2}\int_0^t \sigma(Y_s)d\tilde V_s.
\end{aligned}
\end{equation}

\begin{lemma}\label{lem:scrounge}
\begin{itemize} 
\item[(i)] The stochastic derivatives of the fBm $B^H$ are equal to
\begin{equation*}		
D^{\tilde V}_u B^H _t  = 0, \quad D^{V}_u B^H _t  = K(t,u)\mathbbm 1_{[0,t]}(u).
\end{equation*}
\item[(ii)] The stochastic derivatives of the volatility process $Y$ are
\begin{equation*}
\begin{gathered}
D^{\tilde V}_u Y _t  = 0, 
\\
 D^V_u Y _t  = \left[K(t,u) - \int\limits _u ^t K(s,u) h(s) \exp\left\{ - \int\limits _s ^t h(v)dv \right\} ds \right]\mathbbm 1_{[0,t]}(u)
\end{gathered}
\end{equation*}
where $h(s) := \frac{1}{2} \left(\frac{\kappa}{Y_s ^2}  +\theta \right)$.
\item[(iii)] The stochastic derivatives of $X$ are equal to
\begin{equation*}
\begin{gathered}
D_u^{\tilde V} X _t  = \sqrt{1-\rho^2} \sigma(Y _u) \mathbbm 1_{[0,t]}(u),
\\
D _u^V X _t  =   \left(-  \int_u ^t \sigma(Y _s)\sigma'(Y_s) D _u ^V Y _s ds + \int_u ^t \sigma'(Y_s) D _u ^V Y _s dW_s+\rho \sigma(Y _u)\right)  \mathbbm 1_{[0,t]}(u).
\end{gathered}
\end{equation*}
\end{itemize}
\end{lemma}

Denote
\[
g(y) := f(e^y), \quad F(x) := \int_0^x f(z)dz, \quad G(y) := \int_0^y g(z)dz, \qquad x\ge 0, y\in\mathbb R,
\]
and consider a random variable
\begin{equation}\label{Zoptionprice}
Z_T := \int_0^T \sigma^{-1}(Y_u)d\tilde V_u.
\end{equation}

Note that, due to Assumption \ref{AssumSigma}, (i), $Z_T$ is correctly defined.

\begin{theorem}\label{OptionPriceFormula}
Under Assumptions \ref{AssumSigma} and \ref{AssumPayoff}, the option price $\E f(S_T) = \E g(X_T)$ can be represented as
\begin{equation}\label{OptionPrice1}
\E g(X_T) = \frac{1}{T} \E(G(X_T)Z_T),
\end{equation}
or, alternatively,
\begin{equation}\label{OptionPrice2}
\E f(S_T) = \E\left(\frac{F(S_T)}{S_T}\left(1+\frac{Z_T}{T}\right)\right).
\end{equation}
\end{theorem}


\section{Inverse Euler approximation scheme for the volatility and price processes}\label{sec: Inverse Euler}

Let $0=t^n_0<t^n_1<...<t^n_n = T$ be an equidistant partition of the interval $[0,T]$, $t^n_i = \frac{iT}{n}$, $\Delta_n := \frac{1}{n}$, $\Delta B^H_{k+1}:= B^H_{t^n_{k+1}}- B^H_{t^n_{k}}$ and consider the approximation scheme of the form

\begin{equation}\label{inv. Euler scheme 1}
\hat Y^n_{t^n_{k+1}} = \frac{\hat Y^n_{t^n_{k}} + \frac{\nu}{2}\Delta B^H_{k+1} + \sqrt{(\hat Y^n_{t^n_{k}} + \frac{\nu}{2}\Delta B^H_{k+1})^2 + \kappa\Delta_n(2+\theta\Delta_n)} }{2+\theta\Delta_n}
\end{equation}
with linear interpolation between the points of the partition.

Note that approximations given by \eqref{inv. Euler scheme 1} are strictly positive and it is easy to verify that in points of partition they satisfy the following difference equation:
\begin{equation}\label{diff eq representation}
\hat Y^n_{t^n_{k+1}} 
= \hat Y^n_{t^n_{k}} + \frac{1}{2}\left(\frac{\kappa}{\hat Y^n_{t^n_{k+1}}} - \theta  \hat Y^n_{t^n_{k+1}}\right)\Delta_n + \frac{\nu}{2}\Delta B^H_{k+1}.
\end{equation}

Approximations of the form \eqref{inv. Euler scheme 1} were presented and studied in \cite{HHKW}. We give the result concerning the convergence rate of these approximations (for more detail, see Theorem 4.2 in \cite{HHKW}).

\begin{theorem}\label{China theorem}
Let $\xi\in(0,1)$, $p \ge 2$, $\Delta_n < 1 - \xi$ and parameters $\theta, \kappa, \nu>0$ are such that for all $t\in[0,T]$:
\begin{equation}\label{inv. moment condition}
\kappa\exp\left\{\frac{\theta t}{2}\right\}\ge H(2H-1)(3p+1) \int_0^t \frac{\nu^2}{2}\exp\left\{\frac{\theta s}{2}\right\}|t-s|^{2H-2} ds.
\end{equation}

Then there is such $C=C(T,H, p, Y_0, \theta, \kappa, \nu, \xi)>0$ that
\begin{equation*}
\sup_{t\in[0,T]}\E|Y_t-\hat Y^n_t|^p \le C \Delta_n^{pH}.
\end{equation*}
\end{theorem}

\begin{remark}
Condition \eqref{inv. moment condition} is a sufficient condition for finiteness of the inverse moments of $Y$ of order $3p$, namely for
\begin{equation*}
\sup_{t\in[0,T]}\E \left[\frac{1}{Y_t^{3p}}\right] < \infty.
\end{equation*}
\end{remark}

\begin{figure}[h!]
\centerline{\includegraphics[width=\textwidth]{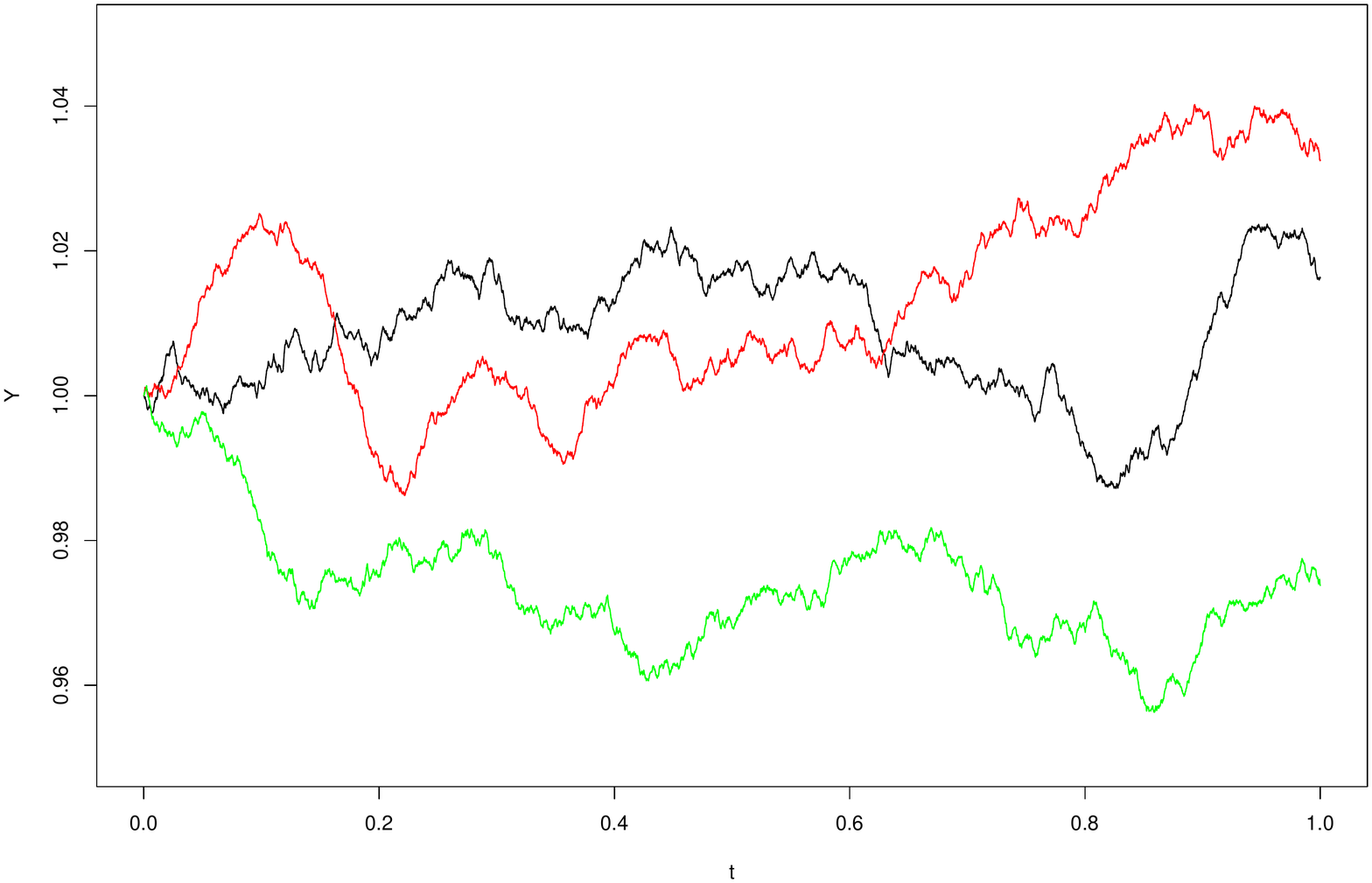}}
\vspace*{8pt}
\caption{Three sample trajectories of the process $Y$ obtained by approximation scheme \eqref{inv. Euler scheme 1}; $T=1$, $\kappa = 1$, $\theta = 1$, $\nu = 0.14$, $Y_0=1$, $H=0.7$ and $\Delta_n = 0.0001$.}\label{Y fol large H}
\end{figure}

Three approximations of the volatility process $Y$ trajectories given by the formula \eqref{inv. Euler scheme 1} with $T=1$, $\kappa = 1$, $\theta = 1$, $\nu = 0.14$, $Y_0=1$, $H=0.7$ and $\Delta_n = 0.0001$ are presented on Fig. \ref{Y fol large H}.

For the sake of simplicity, instead of linear interpolation between the points of the partition, we put $\hat Y^n_t = \hat Y^n_{t^n_{k}}$ for $t\in[\hat Y^n_{t^n_{k}}, \hat Y^n_{t^n_{k+1}})$. It should be noted that in this case speed of convergence of approximations remains the same as in Theorem \ref{China theorem} due to Remark \ref{increments under China condition} because
\begin{equation*}
\E|Y_t-\hat Y^n_{t^n_{k}}|^p \le 2^{p-1} ( \E|Y_{t} - Y_{t^n_{k}}|^p+\E|Y_{t^n_{k}}-\hat Y^n_{t^n_{k}}|^p ).
\end{equation*}

Denote
\begin{equation*}
X_t := \log S_t = X_0 + \mu t  - \frac{1}{2}\int_0^t \sigma^2(Y_s)ds+\int_0^t \sigma(Y_s)dW_s, 
\end{equation*}
where $X_0:=\log S_0$, and consider the discretized process
\begin{equation*}
\begin{aligned}
\hat X_{t_{k}^n}^n &= X_0 + \mu t_{k}^n - \frac{1}{2n}\sum_{j=0}^{k-1} \sigma^2(\hat Y^n_{t_{j}^n})+\sum_{j=0}^{k-1} \sigma(\hat Y^n_{t_{j}^n})\Delta W_j 
\\
&=X_0 + \mu t_{k}^n - \frac{1}{2} \int_0^{t_{k}^n} \sigma^2(\hat Y^n_{s})ds+  \int_0^{t_{k}^n} \sigma(\hat Y^n_{s})dW_s , \quad k = 1,...,n,
\end{aligned}
\end{equation*}
where $\Delta W_j = W_{t_{j+1}^n} - W_{t_{j}^n}$.

Before going to the main theorem of the paper, let us prove several auxiliary results.


\begin{theorem}\label{th moments of approximations}
Let $p\ge 1$. Then, for all $H\in(0,1)$:
\begin{equation*}
\sup_{n\ge 1}\sup_{t\in[0,T]}\E (\hat{Y}^n_t)^p <\infty.
\end{equation*}
\end{theorem}

\begin{figure}[h!]
\centerline{\includegraphics[width=\textwidth]{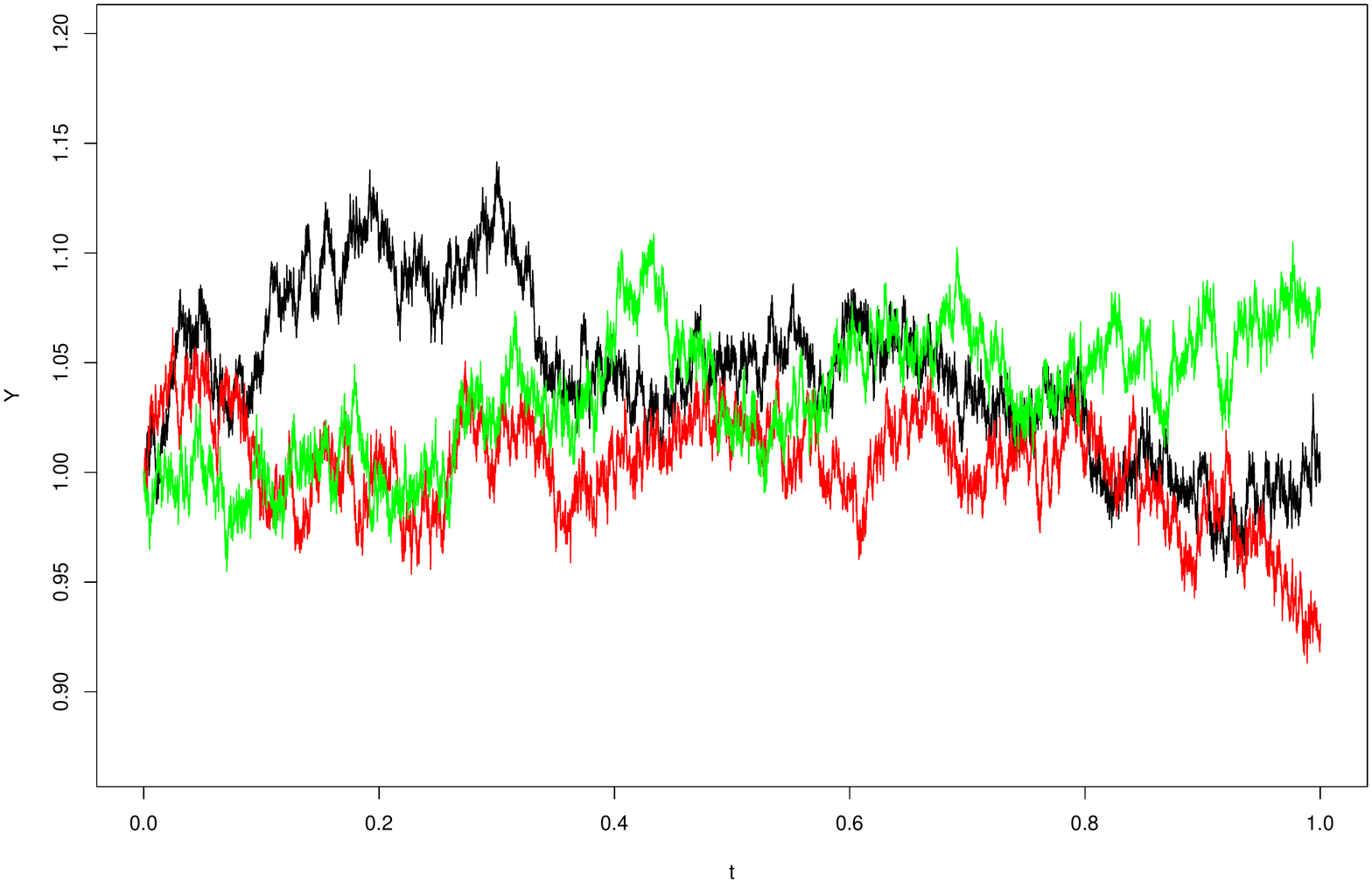}}
\vspace*{8pt}
\caption{Three sample trajectories of the process $\hat Y^n$ for $H=0.3$, $T=1$, $\kappa = 1$, $\theta = 1$, $\nu = 0.14$, $Y_0=1$ and $\Delta_n = 0.0001$.}\label{Y fol small H}
\end{figure}

\begin{remark}
Note that approximations \eqref{inv. Euler scheme 1} (see Fig. \ref{Y fol small H}) are correctly defined for $H<1/2$ and Theorem \ref{th moments of approximations} holds for an arbitrary Hurst parameter as well. However, for $H<1/2$ behaviour of $\hat Y^n$ as $n\to\infty$ remains obscure.
\end{remark}

\begin{corollary}\label{exp moments of approximations}
Approximating processes $\hat Y^n$ have bounded exponential moments, i.e. for any $x>0$ and $\varrho<2$:
\begin{equation*}
\sup_{n\ge1}\E\exp\left\{x\sup_{t\in[0,T]}(\hat Y^n_t)^\varrho\right\}<\infty.
\end{equation*}
\end{corollary}

\begin{remark}\label{moments of approximations}
From Theorem \ref{th moments of approximations}, Corollary \ref{exp moments of approximations} and Assumption \ref{AssumSigma} (ii), using the same argument as in the proof of Theorem \ref{solution to price equation}, it is easy to verify that for any~$m\in\mathbb Z$:
\begin{equation*}
\begin{gathered}
\sup_{n\ge 1} \sup_{t\in[0,T]} \E\left(\hat S^n_t\right)^m < \infty,
\\
\sup_{n\ge 1} \sup_{t\in[0,T]} \E\exp\left\{m \left(\int_0^t \sigma(\hat Y^n_s)dW_s - \frac{1}{2} \int_0^t \sigma^2(\hat Y^n_s)ds \right)\right\} < \infty.
\end{gathered}
\end{equation*}
\end{remark}

\begin{theorem}\label{approx conv theor}
Let $n\ge 2$ and conditions of Theorem \ref{China theorem} hold for $p=4$. Then, under Assumption \ref{AssumSigma}, there exists a constant $C$ such that 
\begin{equation}\label{Xapprox}
\E|X_T - \hat X^n_T|^2 \le C \Delta_n^{2rH},
\end{equation} 
\begin{equation}\label{Zapprox}
\E|Z_T - \hat Z^n_T|^2 \le C \Delta_n^{2rH},
\end{equation}
where $\hat Z^n_T := \int_0^T \sigma^{-1}(\hat Y^n_u) d \tilde V_u$.
\end{theorem}

\begin{lemma}\label{auxiliary approx lemma}
Let $n\ge 2$ and conditions of Theorem \ref{China theorem} hold for $p=32$. Then, under Assumptions \ref{AssumSigma} and \ref{AssumPayoff}, there is such $C_F>0$ that
\[
\E\left|\frac{F(S_T)}{S_T} - \frac{F(\hat S^n_T)}{\hat S^n_T}\right|^2 \le C_F \Delta_n^{H}.
\]
\end{lemma}

\begin{theorem}\label{final option pricing}
Let $n\ge2$ and conditions of Theorem \ref{China theorem} for $p=32$ hold. Then, under Assumptions \ref{AssumSigma} and \ref{AssumPayoff},
\[
\left|\E f(S_T) - \E\left[\frac{F(\hat S^n_T)}{\hat S^n_T}\left(1 + \frac{\hat Z^n_T}{T}\right)\right]\right| \le C\Delta_n^{rH}.
\]
\end{theorem}

\section{Simulations}\label{sec: simulations}

In this section, we use the discretization scheme studied previously to estimate option price for several payoff functions $f$. In all simulations we use $T=1$, $\kappa = 1$, $\theta = 1$ and $\nu = 0.14$ to make sure that for all $H\in(1/2,1)$ the following condition is satisfied for $p=32$:
\[
3p+1 \le \frac{2\kappa}{\nu^2 H T^{2H-1}},
\]
which is sufficient for Theorem \ref{final option pricing} to hold true. For simplicity, we also consider everywhere the case $\mu = 0.5$, $\rho = 0$ and $\sigma = 0.5\left(x+0.01\right)^{0.9}$.

\begin{figure}[h!]
  \centering
\begin{minipage}[b]{0.48\textwidth} \centering
    \includegraphics[width=\textwidth]{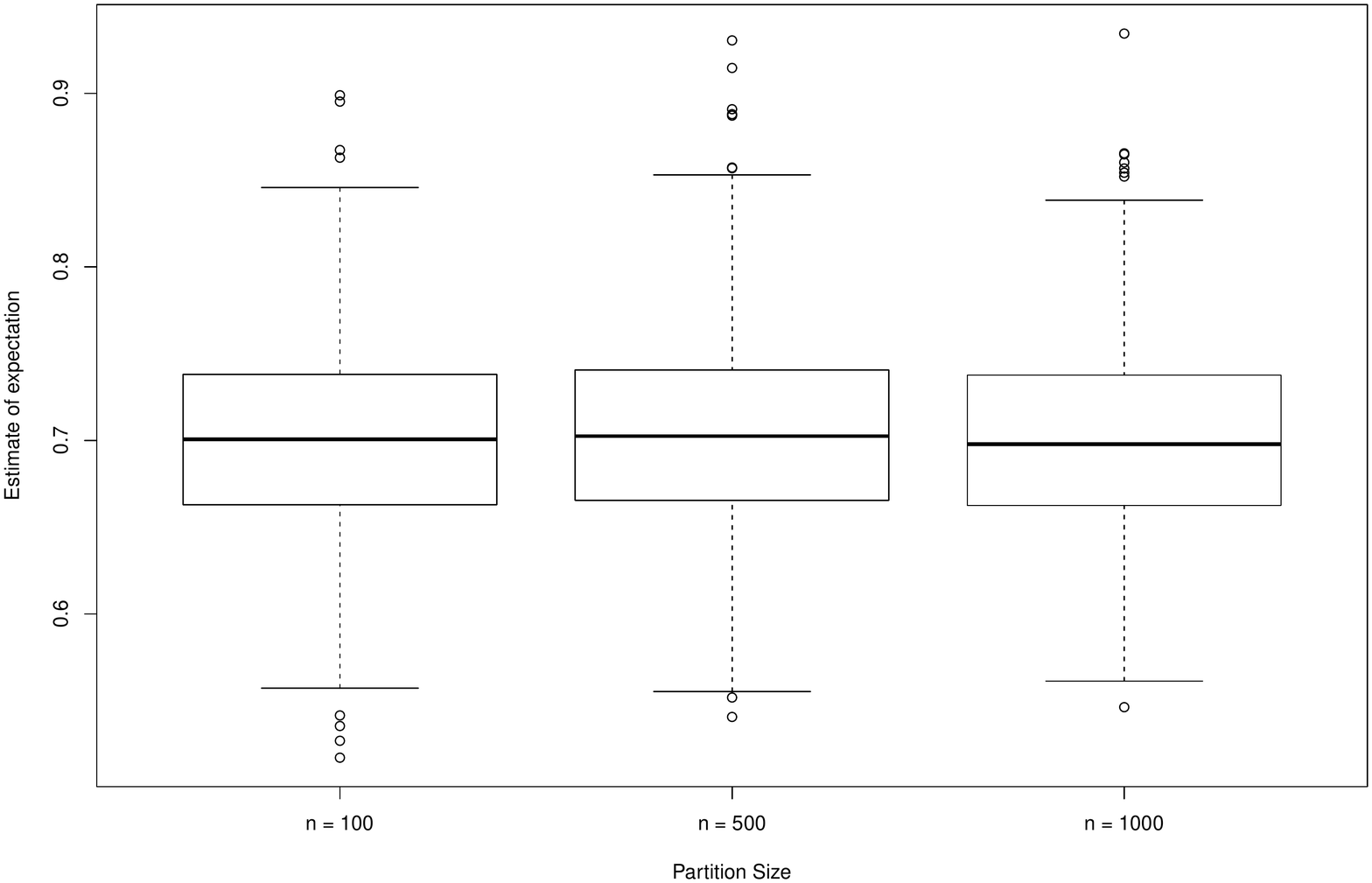}
(a)
\end{minipage}
\begin{minipage}[b]{0.48\textwidth} \centering
\includegraphics[width=\textwidth]{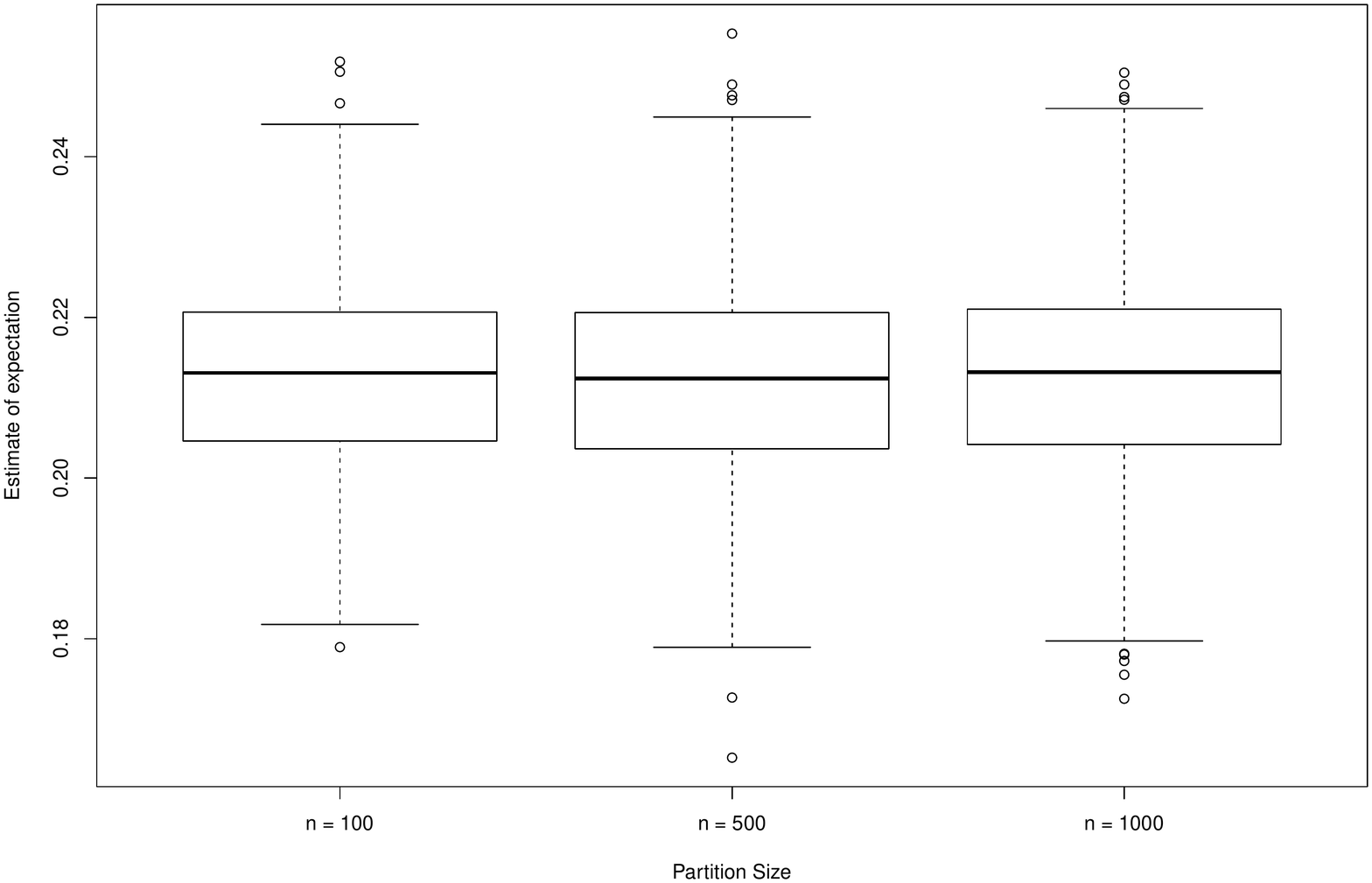}
(b)
\end{minipage}
\begin{minipage}[b]{0.48\textwidth} \centering
\includegraphics[width=\textwidth]{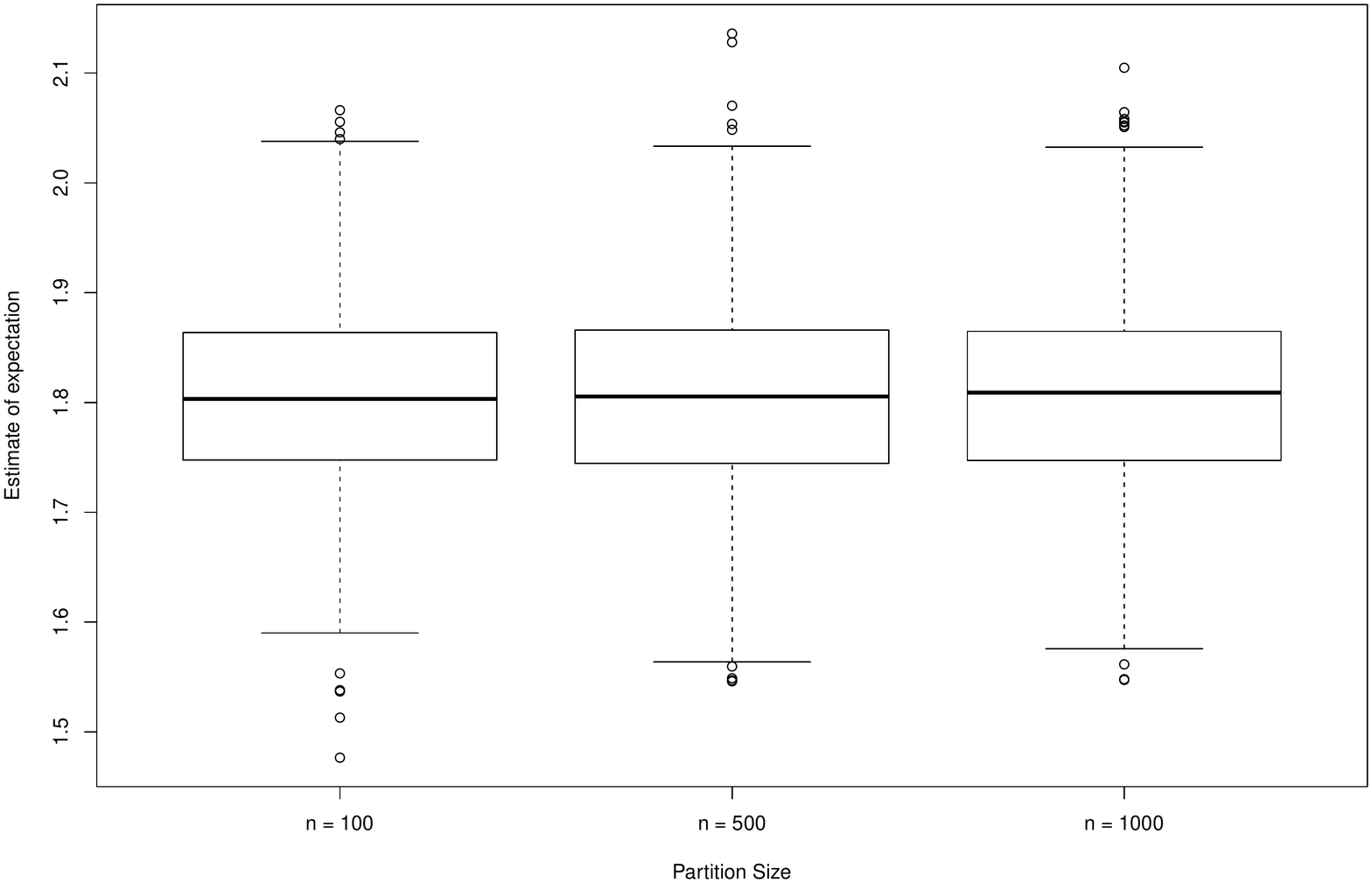}
(c)
\end{minipage}
\caption{Box-and-whisker plots of Monte-Carlo estimates of $\E f(S_T)$ using smoothed formula; in all cases $T=1$, $\kappa = 1$, $\theta = 1$, $\nu = 0.14$, $\mu = 0.5$, $\rho = 0$, $\sigma = 0.5\left(x+0.01\right)^{0.9}$, $H=0.7$; (a) $f(x) = (x - 1)^+$, (b) $f(x) = \mathbbm 1_{[0.5,1]}(x)$, (c) $f(x) = \mathbbm 1_{(0.5,\infty)}(x) + \frac{1}{2} \sum_{k=2}^6 \mathbbm 1_{(0.5k,\infty)}(x)$}\label{Tabl1}
\end{figure}

In Tables 1--3 we present descriptive statistics of Monte-Carlo estimations of $\E\left[\frac{F(\hat S^n_T)}{\hat S^n_T}\left(1 + \frac{\hat Z^n_T}{T}\right)\right]$ (and, therefore, $\E f(S_T)$) for different functions $f$ and different partition sizes $\Delta_n$. On Fig. \ref{Tabl1}, (a)--(c), the data is visualized in a form of box-and-whisker plots. In each case, 1000 Monte-Carlo estimates of option price, calculated from samples of 1000 trials each, were analyzed. All calculations were performed in \textsf{R} using package \textsf{somebm} to generate trajectories of Wiener process and fractional Brownian motion.

\begin{table}[h]
\caption{$f(x) = (x - K)^+$, $K=1$, $\sigma(x) = 0.5(x+0.01)^{0.9}$, $\mu = 0.5$, $H=0.7$}
{\begin{tabular}{ccccccccc} \hline
$n$ & Mean & Standard deviation & Coefficient of variation & Min. & 1st Qu. & Median & 3rd Qu. & Max.\\ \hline\\
100 & 0.7019 & 0.05628101 & 0.0802 & 0.5171 &  0.6630 &  0.7006  &  0.7380 &  0.8989 \\
500 & 0.7040 & 0.05476103 & 0.0778 & 0.5406 &  0.6655 &  0.7025  &  0.7406 &  0.9305\\
1000 & 0.7004 & 0.05459163 & 0.0779 & 0.5463  & 0.6625 &  0.6978  &  0.7375  & 0.9344\\  \hline
\end{tabular}}
\end{table}

\begin{table}[h]
\caption{$f(x) = \mathbbm 1_{[0.5,1]}(x)$, $\sigma(x) = 0.5(x+0.01)^{0.9}$, $\mu = 0.5$, $H=0.7$}
{\begin{tabular}{ccccccccc} \hline
$n$ & Mean & Standard deviation & Coefficient of variation & Min. & 1st Qu. & Median & 3rd Qu. & Max.\\ \hline\\
100 & 0.2126 & 0.01196734 & 0.0563 & 0.1790 & 0.2046 &  0.2131  &  0.2206 &  0.2518 \\
500 & 0.2123 & 0.01266216 & 0.0596 & 0.1652 &  0.2037 &  0.2124 &  0.2206 &  0.2553\\
1000 & 0.2129 & 0.01272749 & 0.0598 & 0.1725 & 0.2042  &  0.2132  &  0.2210 &  0.2505\\  \hline
\end{tabular}}
\end{table}

\begin{table}[h]
\caption{$f(x) = \mathbbm 1_{(0.5,\infty)}(x) + \frac{1}{2} \sum_{k=2}^6 \mathbbm 1_{(0.5k,\infty)}(x)$, $\sigma(x) = 0.5(x+0.01)^{0.9}$, $\mu = 0.5$, $H=0.7$}
{\begin{tabular}{ccccccccc} \hline
$n$ & Mean & Standard deviation & Coefficient of variation & Min. & 1st Qu. & Median & 3rd Qu. & Max.\\ \hline\\
100 & 1.804 & 0.08973507 & 0.0497 & 1.476 &  1.748  &  1.803   &   1.864  &  2.066 \\
500 & 1.806 & 0.08873267 & 0.0491 & 1.546 & 1.745  & 1.805   &   1.866  &  2.136\\
1000 & 1.806 & 0.09001699 & 0.0498 & 1.547 &  1.747  & 1.809   &   1.865 &  2.105\\  \hline
\end{tabular}}
\end{table}

As we can see, simulations show relatively small coefficient of variation in all cases. Nota that increasing partition size does not lead to any significant changes in standard deviation of the estimates.

\section{Proofs}\label{sec: proofs}

\textbf{Proof of Theorem \ref{th. int. inv. m.}.} Denote $\alpha:=1-\beta$ and let $\varepsilon>0$ be fixed. By applying the chain rule, we obtain:
\begin{equation}\label{inv. m: chain rule}
\begin{aligned}
(Y_t+\varepsilon)^\alpha &= (Y_{t_0}+\varepsilon)^\alpha + \int_{t_0}^t \frac{\kappa\alpha}{2Y_u(Y_u+\varepsilon)^{1-\alpha}}du - \int_{t_0}^t \frac{\theta\alpha Y_u}{2(Y_u+\varepsilon)^{1-\alpha}}du
\\
&\quad + \int_{t_0}^t \frac{\nu \alpha}{2(Y_u+\varepsilon)^{1-\alpha}}dB_u^H.
\end{aligned}
\end{equation}

It is clear from \eqref{fHestonVolatility} that the process $Y = \{Y_t,~t\in[0,T]\}$ has trajectories that are $\delta$-H\"older-continuous for any $\delta\in(0,H)$, so the process 
\begin{equation*}
\frac{\nu\alpha}{2(Y_t+\varepsilon)^{1-\alpha}},\quad t\in[0,T], 
\end{equation*}
also has H\"older-continuous trajectories up to the order $H$. Therefore, the sum of H\"older exponents of the integrator and integrand in the integral w.r.t. fractional Brownian motion in \eqref{inv. m: chain rule} exceeds 1. In this case this integral is the pathwise limit of Riemann-Stieltjes integral sums (see, for example, \cite{Zahle}), coincides with the pathwise Stratonovich integral and, by applying Theorem \ref{NualTheorem}, we can rewrite \eqref{inv. m: chain rule} as follows:
\begin{equation}\label{inv. m: Nual}
\begin{aligned}
(Y_t+\varepsilon)^\alpha &= (Y_{t_0}+\varepsilon)^\alpha+\int_{t_0}^t \frac{\kappa\alpha}{2Y_u(Y_u+\varepsilon)^{1-\alpha}}ds - \int_{t_0}^t \frac{\theta\alpha Y_u}{2(Y_u+\varepsilon)^{1-\alpha}}ds
\\
&\quad+ H(2H-1)\int_{t_0}^t\int_0^t D_s^H\left[\frac{\nu\alpha}{2(Y_u+\varepsilon)^{1-\alpha}}\right]|u-s|^{2H-2}ds du
\\
&\quad + \int_{t_0}^t \frac{\nu \alpha}{2(Y_u+\varepsilon)^{1-\alpha}}\delta B_u^H,
\end{aligned}
\end{equation}
where $D^H_s$ is the Malliavin derivative operator w.r.t. $B^H$ and $\int_{t_0}^t \frac{\nu \alpha}{2(Y(s)+\varepsilon)^{1-\alpha}}\delta B_s^H$ is the corresponding Skorokhod integral.

Note that
\begin{equation*}
\begin{aligned}
D_s^H Y_u = & D_s^H\left[Y_0+\int_0^u \frac{\kappa}{2Y_v}dv - \frac{\theta}{2}\int_0^u Y_v dv +\frac{\nu}{2}B_u^H\right] = 
\\
=&-\int_0^u \frac{\kappa D_s^H Y_v }{2Y^2_v}dv - \frac{\theta}{2}\int_0^u D_s^H Y_v dv +\frac{\nu}{2}\mathbbm 1_{[0,u]}(s)
\\
=& - \int_0^u \left(\frac{\kappa}{2Y^2_v}+\frac{\theta}{2}\right)D_s^H Y_v dv + \frac{\nu}{2}\mathbbm 1_{[0,u]}(s).
\end{aligned} 
\end{equation*}

From this, it is easy to verify that
\begin{equation*}
D_s^H Y_u = \frac{\nu}{2}\exp\left\{-\int_s^u \left(\frac{\kappa}{2Y^2_v}+\frac{\theta}{2}\right)dv \right\}\mathbbm 1_{[0,u]}(s),
\end{equation*}
so
\begin{equation}\label{inv. m: Mal. der}
\begin{aligned}
D_s^H \left[\frac{\nu\alpha}{2(Y_u+\varepsilon)^{1-\alpha}}\right] &= -\frac{\nu\alpha(1-\alpha)}{2(Y_u+\varepsilon)^{2-\alpha}} D_s^H Y_u
\\
&= -\frac{\nu^2\alpha(1-\alpha)}{4(Y_u+\varepsilon)^{2-\alpha}}\exp\left\{-\int_s^u \left(\frac{\kappa}{2Y^2_v}+\frac{\theta}{2}\right)dv \right\}\mathbbm 1_{[0,u]}(s).
\end{aligned} 
\end{equation}

Taking into account \eqref{inv. m: Nual} and \eqref{inv. m: Mal. der}, we can rewrite \eqref{inv. m: chain rule} in the following form:
\begin{equation}\label{inv. m: to bound}
\begin{aligned}
(Y_t+\varepsilon)^\alpha &= (Y_{t_0}+\varepsilon)^\alpha+\int_{t_0}^t \frac{\kappa\alpha}{2Y_u(Y_u+\varepsilon)^{1-\alpha}}ds - \int_{t_0}^t \frac{\theta\alpha Y_u}{2(Y_u+\varepsilon)^{1-\alpha}}ds
\\
&\quad- \int_{t_0}^t\int_0^u \frac{\frac{\nu^2}{2}\alpha(1-\alpha)\exp\left\{-\int_s^u \left(\frac{\kappa}{2Y^2_v}+\frac{\theta}{2}\right)dv \right\} \varphi(u,s) ds }{2(Y_u+\varepsilon)^{2-\alpha}} du
\\
&\quad + \int_{t_0}^t \frac{\nu \alpha}{2(Y_u+\varepsilon)^{1-\alpha}}\delta B_u^H,
\end{aligned}
\end{equation}
where $\varphi(u,s):=H(2H-1)|u-s|^{2H-2}$.

Note that 
\begin{equation}\label{inv. m: ub1}
\begin{gathered}
\int_{t_0}^t \frac{\kappa\alpha}{2Y_u(Y_u+\varepsilon)^{1-\alpha}}ds - \int_{t_0}^t\int_0^u \frac{\frac{\nu^2}{2}\alpha(1-\alpha)\exp\left\{-\int_s^u \left(\frac{\kappa}{2Y^2_v}+\frac{\theta}{2}\right)dv \right\} \varphi(u,s) ds }{2(Y_u+\varepsilon)^{2-\alpha}} du
\\
\ge \alpha \int_{t_0}^t \frac{\kappa-(1-\alpha)\int_0^u \frac{\nu^2}{2} \exp\left\{-\int_s^u \left(\frac{\kappa}{2Y^2_v}+\frac{\theta}{2}\right)dv \right\} \varphi(u,s) ds }{2(Y_u+\varepsilon)^{2-\alpha}}du.
\end{gathered}
\end{equation}

It is easy to verify that
\begin{equation*}
\begin{aligned}
0 &\ge -\int_0^u \exp\left\{-\int_s^u \left(\frac{\kappa}{2Y^2_v}+\frac{\theta}{2}\right)dv \right\} \varphi(u,s) ds
\\
&\ge - H(2H-1) \int_0^u |u-s|^{2H-2} ds
\\
&\ge - H T^{2H-1},
\end{aligned}
\end{equation*}
so
\begin{equation*}
0 \ge -(1-\alpha)\int_0^u \frac{\nu^2}{2} \exp\left\{-\int_s^u \left(\frac{\kappa}{2Y^2_v}+\frac{\theta}{2}\right)dv \right\} \varphi(u,s) ds \ge - (1-\alpha)\frac{\nu^2}{2} H T^{2H-1}.
\end{equation*}

Hence, if $\alpha\in\left(\max\{0, 1- \frac{\kappa}{\nu^2 H T^{2H-1}}\}, 1\right)$, i.e. when $0\ge - (1-\alpha)\frac{\nu^2}{2} H T^{2H-1}\ge - \frac{\kappa}{2}$,
\begin{equation*}
\kappa-(1-\alpha)\int_0^u \frac{\nu^2}{2} \exp\left\{-\int_s^u \left(\frac{\kappa}{2Y^2_v}+\frac{\theta}{2}\right)dv \right\} \varphi(u,s) ds \ge \frac{\kappa}{2},
\end{equation*}
and
\begin{equation}\label{inv. m: ub2}
\begin{gathered}
\int_{t_0}^t \frac{\kappa\alpha}{2Y_u(Y_u+\varepsilon)^{1-\alpha}}ds - \int_{t_0}^t\int_0^u \frac{\frac{\nu^2}{2}\alpha(1-\alpha)\exp\left\{-\int_s^u \left(\frac{\kappa}{2Y^2_v}+\frac{\theta}{2}\right)dv \right\} \varphi(u,s) ds }{2(Y_u+\varepsilon)^{2-\alpha}} du
\\
\ge \frac{\alpha \kappa}{4} \int_{t_0}^t \frac{1}{(Y_u+\varepsilon)^{2-\alpha}}du.
\end{gathered}
\end{equation}

Moreover,
\begin{equation}\label{inv. m: ub3}
-\int_{t_0}^t \frac{\theta\alpha Y_u}{2(Y_u+\varepsilon)^{1-\alpha}}ds \ge - \frac{\theta\alpha}{2} \int_{t_0}^t (Y_u+\varepsilon)^\alpha du. 
\end{equation}

Therefore, taking into account upper bounds \eqref{inv. m: ub1}, \eqref{inv. m: ub2} and \eqref{inv. m: ub3}, it is obvious from \eqref{inv. m: to bound} that
\begin{equation*}
\begin{aligned}
(Y_t+\varepsilon)^\alpha &\ge (Y_{t_0}+\varepsilon)^\alpha + \frac{\alpha \kappa}{4} \int_{t_0}^t \frac{1}{(Y_u+\varepsilon)^{2-\alpha}}du - \frac{\theta\alpha}{2} \int_{t_0}^t (Y_u+\varepsilon)^\alpha du
\\
&\quad + \int_{t_0}^t \frac{\nu \alpha}{2(Y_u+\varepsilon)^{1-\alpha}}\delta B_u^H,
\end{aligned}
\end{equation*}
or
\begin{equation*}
\begin{aligned}
\int_{t_0}^t \frac{1}{(Y_u+\varepsilon)^{2-\alpha}}du \le& \frac{4}{\kappa\alpha}(Y_t+\varepsilon)^\alpha - \frac{4}{\kappa\alpha}(Y_{t_0}+\varepsilon)^\alpha + \frac{2\theta}{\kappa} \int_{t_0}^t (Y_u+\varepsilon)^\alpha du 
\\
& -  \frac{4}{\kappa\alpha} \int_{t_0}^t \frac{\nu \alpha}{2(Y_u+\varepsilon)^{1-\alpha}}\delta B_u^H.
\end{aligned}
\end{equation*}

Since the expectation of the Skorokhod integral is zero, by letting $\varepsilon\to 0$ we obtain that
\begin{equation}\label{upper bound for IIM}
\E \left[\int_{t_0}^t \frac{1}{Y_u^{2-\alpha}}du\right] \le \frac{4}{\kappa\alpha}\E (Y_t^\alpha -  Y_{t_0}^\alpha )  + \frac{2\theta}{\kappa} \int_{t_0}^t \E Y_u^\alpha du.
\end{equation}

Finiteness of the right-hand side of \eqref{upper bound for IIM} follows from Theorem \ref{UpperBoundForY}. \QED
\\[11pt]
\textbf{Proof of Theorem \ref{theor increments}.}  From \eqref{fHestonVolatility}, H\"older's and Jensen's inequalities it is clear that
\begin{equation}\label{incr 1}
\begin{gathered}
\mathbb E |Y_t - Y_s|^{1+\beta} = \E\left|\int_s^t \frac{\kappa}{2Y_u} du - \frac{\theta}{2} \int_s^t Y_u du +\frac{\nu}{2}(B^H_t - B^H_s) \right|^{1+\beta}
\\
\le 3^\beta \left(\left(\frac{\kappa}{2}\right)^{1+\beta}\E\left|\int_s^t \frac{1}{Y_u} du\right|^{1+\beta} + \left(\frac{\theta}{2}\right)^{1+\beta}\E\left|\int_s^t Y_u du\right|^{1+\beta} + \left(\frac{\nu}{2}\right)^{1+\beta}\E\left|B^H_t - B^H_s\right|^{1+\beta}\right)
\\
\le \tilde C_1 (t-s)^\beta \E\int_s^t \frac{1}{Y_u^{1+\beta}}du + \tilde C_2 (t-s)^\beta \E \int_s^t Y_u^{1+\beta} du + C_{3} |t-s|^{(1+\beta)H},
\end{gathered}
\end{equation}
where
\begin{equation*}
\tilde C_1 := 3^\beta\left(\frac{\kappa}{2}\right)^{1+\beta}, \quad \tilde C_2 := 3^\beta\left(\frac{\theta}{2}\right)^{1+\beta},\quad
C_3 := (3\sqrt{2})^\beta\sqrt{\frac{2}{\pi}}\Gamma\left(1+\frac{\beta}{2}\right)\left(\frac{\nu}{2}\right)^{1+\beta}.
\end{equation*}
Note that form of $C_3$ follows from the fact that $B^H_t - B^H_s \sim \mathcal N(0, |t-s|^{2H})$ (see, for example, \cite{Winkel}) 

From Theorem \ref{UpperBoundForY} it is obvious that
\begin{equation}\label{incr 2}
 \tilde C_2 (t-s)^\beta \E \int_s^t Y_u^{1+\beta} du \le \tilde C_2 \sup_{u\in[0,T]}\E\left[ Y_u^{1+\beta}\right] (t-s)^{1+\beta} =: C_2 (t-s)^{1+\beta}.
\end{equation}
 
Finally, from Theorem \ref{th. int. inv. m.},
\begin{equation}\label{incr 3}
\tilde C_1 (t-s)^\beta \E\int_s^t \frac{1}{Y_u^{1+\beta}}du < C_1 (t-s)^\beta,
\end{equation}
where $C_1 = \tilde C_1 \E\int_s^t \frac{1}{Y_u^{1+\beta}}du < C_1 (t-s)^\beta$.

The statement of the Theorem now follows from \eqref{incr 1}, \eqref{incr 2} and \eqref{incr 3} as well as the fact that from condition $\beta<1$ it is easy to verify that for any $H\in(1/2, 1)$:
\begin{equation*}
\beta < (1+\beta)H .
\end{equation*} \QED
\\[11pt]
\textbf{Proof of Theorem \ref{solution to price equation}.} 1. From Theorem \ref{UpperBoundForY}, for all $\varrho\in[0,2)$:
\begin{equation*}
\sup_{t\in[0,T]}|Y(t)|^\varrho \le C_1 + C_2 \sup_{t\in [0,T]} |B^H_t|^\varrho,
\end{equation*}
and, due to \cite{Fernique}, for all $x>0$ and $\varrho\in[0,2)$:
\begin{equation*}
\mathbb E \exp\left\{x\sup_{t\in [0,T]} |B^H_t|^\varrho\right\}<\infty.
\end{equation*}

Hence,
\begin{equation*}
\exp\left\{x\sup_{t\in [0,T]} |Y_t|^\varrho\right\} \le \exp\left\{C_1x + C_2 x\sup_{t\in [0,T]} |B^H_t|^\varrho\right\} < \infty.
\end{equation*}

2. In order to show that the representation \eqref{ExplicitRepresentationOfS} indeed holds, it is sufficient to prove that the integrals $\int_0^t \sigma(Y_s)dW_s$ and $\int_0^t \sigma(Y_s)S_s dW_s$ are well-defined, while the form of the representation can be obtained straightforwardly.

Note that (see, for example, \cite{MishuraFBM}) for all $p>0$
\begin{equation*}
\mathbb{E} \sup_{s\in [0,T]} |B^H_s|^{p} <\infty,
\end{equation*}
so, due to item $(ii)$ from Assumption \ref{AssumSigma} and Theorem \ref{UpperBoundForY},
\begin{equation*}
\int_0^t \mathbb{E} \sigma^2(Y_s)ds \le C^2_\sigma \int_0^t \mathbb{E} (1+|Y_s|^q)^2ds \le 2 C^2_\sigma \int_0^t \mathbb{E} (1+|Y_s|^{2q})ds < \infty,
\end{equation*}
and the integral $\int_0^t \sigma(Y_s)dW_s$ is well-defined.

Now consider the integral $\int_0^t \sigma(Y_s)S_s dW_s$. As
\begin{equation}\label{IntSigmaS}
\begin{gathered}
\int_0^T \mathbb E\left[\sigma^2(Y_s)S^2_s\right]ds \le \int_0^T\left( \mathbb E \sigma^4(Y_s)\right)^\frac{1}{2}\left( \mathbb E S^4_s\right)^\frac{1}{2}ds
\\
\le T\sup_{s\in[0,T]} \left( \mathbb E \sigma^4(Y_s)\right)^\frac{1}{2} \sup_{s\in[0,T]} \left( \mathbb E S^4_s\right)^\frac{1}{2},
\end{gathered}
\end{equation}
it is sufficient to check two conditions:
\begin{equation*}
\begin{gathered}
\sup_{s\in[0,T]} \left( \mathbb E \sigma^4(Y_s)\right)^\frac{1}{2} <\infty,\qquad \sup_{s\in[0,T]} \left( \mathbb E S^4_s\right)^\frac{1}{2} <\infty.
\end{gathered}
\end{equation*}

Using Theorem \ref{UpperBoundForY} and Assumption \ref{AssumSigma}, $(ii)$, it is easy to verify that
\begin{equation}\label{UpperBoundFor4}
\begin{gathered}
\sup_{s\in[0,T]}\left(\mathbb E \sigma^4 (Y_s)\right)^\frac{1}{2} \le \sup_{s\in[0,T]}\left( C_\sigma^{4}\mathbb E (1+|Y_s|^q)^{4}\right)^\frac{1}{2} 
\\
\le \sup_{s\in[0,T]} \left( 8 C_\sigma^{4}\mathbb E (1+|Y_s|^{4q}) \right)^\frac{1}{2} \le \left(8 C_\sigma^4 +  8 C_\sigma^4\sup_{s\in[0,T]}\mathbb E |Y_s|^{4q}\right)^\frac{1}{2}<\infty.
\end{gathered}
\end{equation} 

Moreover, from \eqref{YExponentialMoments}, for any $x>0$:
\begin{equation}\label{GenerNovikov}
\begin{aligned}
\mathbb E \exp\left\{x \int_0^t \sigma^2(Y_s)ds \right\} &\le \mathbb{E} \exp\left\{2xC_\sigma   \int_0^t(1+|Y_s|^{2q})ds\right\}
\\
&\le C\mathbb E \exp\left\{2xC_\sigma T \sup_{s\in[0,T]} |Y_s|^{2q}\right\} <\infty,
\end{aligned}
\end{equation}
hence, for all $n\in\mathbb{Z}$, by putting $x:=\frac{4n^2}{2}$, we obtain the Novikov's condition for the process $-2n\sigma(Y_t)$, $t\in~[0,T]$.

Consequently,
\begin{equation}\label{DDexpY}
\mathbb E \exp\left\{2n\int_0^t \sigma(Y_s)dW_s - 2n^2 \int_0^t \sigma^2(Y_s)ds\right\} = 1,
\end{equation}
and so
\begin{equation}\label{UpperBoundForSn}
\begin{gathered}
\sup_{t\in[0,T]} \mathbb E S_t^n \le C\sup_{t\in[0,T]}\mathbb E\exp\left\{n\int_0^t \sigma(Y_s)dW_s - \frac{n}{2}\int_0^t\sigma^2(Y_s)ds\right\}
\\
= C\sup_{t\in[0,T]}\mathbb E\left[\exp\left\{n\int_0^t \sigma(Y_s)dW_s - n^2 \int_0^t\sigma^2(Y_s)ds\right\}\exp\left\{(n^2-\frac{n}{2})\int_0^t\sigma^2(Y_s)ds\right\}\right]
\\
\le C\sup_{t\in[0,T]} \Bigg[\left(\mathbb E\exp\left\{2n\int_0^t \sigma(Y_s)dW_s - n\int_0^t\sigma^2(Y_s)ds\right\}\right)^{\frac{1}{2}}\times
\\
\times \left(\mathbb E\exp\left\{(2n^2-n)\int_0^t\sigma^2(Y_s)ds\right\}\right)^{\frac{1}{2}}\Bigg]
\\
= C \sup_{t\in[0,T]} \left(\mathbb E\exp\left\{(2n^2-n)\int_0^t\sigma^2(Y_s)ds\right\}\right)^{\frac{1}{2}} <\infty
\end{gathered}
\end{equation}
due to \eqref{GenerNovikov}.

Therefore, from \eqref{IntSigmaS}, \eqref{UpperBoundFor4} and \eqref{UpperBoundForSn},
\begin{equation*}
\int_0^T \mathbb E\left[\sigma^2(Y_s)S^2_s\right]ds < \infty
\end{equation*}
and so the integral $\int_0^t \sigma(Y_s)S_s dW_s$ is well-defined. \QED
\\[11pt]
\textbf{Proof of Theorem \ref{th. arb-free}.} The proof is similar to the proof of Theorem 4 in \cite{BdPM}.
\\[11pt]
\textbf{Proof of Lemma \ref{lem:scrounge}.} Item $(i)$ can be found in \cite{BdPM}. In particular, $D^{\tilde V}_u Y _t  = 0$ in $(ii)$ follows from independence of $Y$ and $V$.

Applying stochastic derivative operator to both parts of the integral form of \eqref{fHestonVolatility}, we get
\begin{equation}\label{cognoscenti}
\begin{aligned}
D^V_u Y_t &= \frac 12 \int _{0} ^t D ^V _u \left(\frac{\kappa}{Y_s} - \theta Y _s\right)ds + \nu D ^V _u  B ^ H _t
\\
&=
- \frac 12 \int\limits _{0} ^t   \left(\frac{\kappa}{Y ^2 _s} + \theta \right) D ^V _u {Y_s}ds 
+ \nu K(t,u)\mathbbm 1_{[0,t]}(u)
\\
& = - \frac 12 \int\limits _{0} ^t   h(s) D ^V _u {Y_s}ds + \nu K(t,u)\mathbbm 1_{[0,t]}(u).
\end{aligned}
\end{equation}

Application of the chain rule with the function $F(x) = 1/x$  can be justified by the same argument as in Remark 10 of \cite{BdPM}, since $F$ is locally Lipschitz on $(0,\infty)$.

According to \cite{MYuT1}, Theorem 2, $Y$ does not hit zero a.s. Therefore $h$ is well defined a.s., and  \eqref{cognoscenti}
means that for a fixed $u$, the process $\{Z_t, t \in [0,t]\}$
defined by $Z _t :=  D ^B _u Y_t$ satisfies a random linear integral equation of the form
\begin{equation}\label{smug}
Z_t = -  \int\limits _{0} ^t h(s) Z_s  ds + \nu K(t,u)\mathbbm 1_{[0,t]}(u).
\end{equation}

This is a Volterra equation, and its solution is given by 
\begin{equation}\label{profligacy}
Z _t = \nu \left[K(t,u) - \int\limits _u ^t K(s,u) h(s)  \exp \left\{ - \int\limits _{s} ^t h(v)dv\right\} ds \right]\mathbbm 1_{[0,t]}(u).
\end{equation}

Note that $K$ is differentiable in the first argument ($\frac{\partial}{dt}K(t,s)$ is well defined for $t>s$), so \eqref{profligacy} can be checked by 
substituting in \eqref{smug} and taking derivatives of both sides.
	
Both derivatives in $(iii)$ are obtained by direct differentiation following the Malliavin derivative rules, see e.g. \cite{NP88}, Proposition 3.4. Since $Y$ is independent of $\tilde V$,
	
	\[
	D _u^{\tilde V} X _t  =  \sqrt{1-\rho^2} D _u^{\tilde V}  \int_0 ^t \sigma(Y _s) d\tilde V_s =  \sqrt{1-\rho^2} \sigma(Y _u) \mathbbm 1_{[0,t]}(u).
	\]

To find $D _u^V X _t$, we note that
	
\begin{equation*}
\begin{aligned}
D _u^V X _t  &=   D _u^V \left[- \frac{1}{2}\int_0^t \sigma^2(Y_s)ds + \sqrt{1-\rho^2}\int_0^t \sigma(Y_s)d\tilde V_s + \rho\int_0^t \sigma(Y_s)dV_s\right]
\\
		&=  - \int_0 ^t \sigma(Y _s)\sigma'(Y_s) D _u^V Y _s  ds + \sqrt{1-\rho^2} \int_0 ^t \sigma'(Y_s) D _u^V Y _s d\tilde V_s
\\
		&\quad+ \rho \int_0^t \sigma'(Y_s) D _u^V Y _s dV_s + \rho \sigma(Y_u) \mathbbm 1_{[0,t]}(u)
\\
		&= \left(-  \int_u ^t \sigma(Y _s)\sigma'(Y_s) D _u ^V Y _s ds + \int_u ^t \sigma'(Y_s) D _u ^V Y _s dW_s+\rho \sigma(Y _u)\right)  \mathbbm 1_{[0,t]}(u).
\end{aligned}
\end{equation*}
\QED
\\[11pt]
\textbf{Proof of Theorem \ref{OptionPriceFormula}.}  The result can be obtained by following the proof of Lemma 11 in \cite{BdPM}, taking into account Lemma \ref{lem:scrounge} and relation \eqref{UpperBoundForSn}. 
\QED
\\[11pt]
\textbf{Proof of Theorem \ref{th moments of approximations}.} First, note that for any fixed $n$ and $k=0,1,...,n$:
\begin{equation}\label{finite n moments}
\begin{gathered}
(\hat Y^n_{t^n_{k+1}})^p = \left(\frac{\hat Y^n_{t^n_{k}} + \frac{\nu}{2}\Delta B^H_{k+1} + \sqrt{(\hat Y^n_{t^n_{k}} + \frac{\nu}{2}\Delta B^H_{k+1})^2 + \kappa\Delta_n(2+\theta\Delta_n)} }{2+\theta\Delta_n}\right)^p
\\
\le C \left((\hat Y^n_{t^n_{k}})^p + \left(\frac{\nu}{2}\right)^p |\Delta B^H_{k+1}|^p + \left((\hat Y^n_{t^n_{k}} + \frac{\nu}{2}\Delta B^H_{k+1})^2 + \kappa\Delta_n(2+\theta\Delta_n)\right)^{\frac{p}{2}}\right)
\\
\le C \left((\hat Y^n_{t^n_{k}})^p + |\Delta B^H_{k+1}|^p + |\hat Y^n_{t^n_{k}} + \frac{\nu}{2}\Delta B^H_{k+1}|^p + (\kappa\Delta_n(2+\theta\Delta_n))^{\frac{p}{2}}\right)
\\
\le C\left(1+(\hat Y^n_{t^n_{k}})^p + |\Delta B^H_{k+1}|^p\right)
\\
\le C\left(1+(\hat Y^n_{t^n_{k}})^p + \sup_{t\in[0,T]}|B^H_{t}|^p\right).
\end{gathered}
\end{equation}

By continuing calculations above recurrently and taking into account that $\hat Y^n_{t^n_{0}} = Y_0$, it is easy to see that there is such constant $C_n$ that
\[
\sup_{t\in[0,T]}(\hat Y^n_t)^p = \max_{k=0,...,n} (\hat Y^n_{t^n_{k}})^p < C_n(1+  \sup_{t\in[0,T]}|B^H_{t}|^p).
\]

Moreover, for any fixed $N$ there is such constant $C_N$ that 
\[
\sup_{1\le n\le N}\sup_{t\in[0,T]}\E (\hat{Y}^n_t)^p = \max_{n=1,...,N}\max_{k=0,...,n} \E(\hat Y^n_{t^n_{k}})^p < C_N(1+  \sup_{t\in[0,T]}|B^H_{t}|^p).
\]

Let us prove that there is such $C>0$ (which does not depend on $n$) that
\[
\sup_{n\ge 1}\sup_{t\in[0,T]}(\hat{Y}^n_t)^p < C(1+  \sup_{t\in[0,T]}|B^H_{t}|^p).
\]

From calculations above, it will be enough to show that, for some $N\ge 1$,
\[
\sup_{n > N}\sup_{t\in[0,T]}(\hat{Y}^n_t)^p < C(1+  \sup_{t\in[0,T]}|B^H_{t}|^p).
\]

Let $n>2(8\theta)^pT^{p-1}$ be fixed. Consider the last moment of staying above level $Y_0/2$, i.e. 
\begin{equation*}
\tau_1 := \max\left\{k=1,...,n~|~\forall t^n_l \le t^n_k: \hat Y^n_{t^n_l} \ge \frac{Y_0}{2}\right\}.
\end{equation*}

Let us prove that for any point of the partition $t^n_k$, $k=1,...,n$, the following inequality holds:
\begin{equation}\label{pre-Gronwall}
\begin{aligned}
(\hat Y^n_{t^n_k})^p \le& \left((4Y_0)^p + \left(\frac{8\kappa T}{Y_0}\right)^p + (8\nu)^p \sup_{s\in[0,T]}|B^H_s|^p\right) 
\\
&\quad+ (8\theta)^p T^{p-1}\sum_{j=1}^{k} (\hat Y^n_{t^n_j})^p \Delta_n.
\end{aligned}
\end{equation}

In order to do that, we will separately consider cases $t^n_k \le t^n_{\tau_1}$ and $t^n_k > t^n_{\tau_1}$.

\textbf{Step 1.} Assume that $t^n_k \le t^n_{\tau_1}$. Then, due to representation \eqref{diff eq representation},
\begin{equation*}
\begin{gathered}
(\hat Y^n_{t^n_k})^p = \left(Y_0 + \frac{1}{2}\sum_{j=1}^{k}\left(\frac{\kappa}{\hat Y^n_{t^n_j}} - \theta \hat Y^n_{t^n_j}\right) \Delta_n +\frac{\nu}{2} B^H_{t^n_k}\right)^p
\\
\le 4^{p-1} \left(Y_0^p + \left(\frac{1}{2}\sum_{j=1}^{k} \frac{\kappa}{\hat Y^n_{t^n_j}} \Delta_n\right)^p + \left(\frac{\theta}{2}\sum_{j=1}^{k} \hat Y^n_{t^n_j} \Delta_n\right)^p + \left(\frac{\nu}{2}\right)^p |B^H_{t^n_k}|^p\right).
\end{gathered}
\end{equation*}

Note that for all $t^n_k \le t^n_{\tau_1}$:
\begin{equation*}
\left(\frac{1}{2}\sum_{j=1}^{k} \frac{\kappa}{\hat Y^n_{t^n_j}} \Delta_n\right)^p \le \left(\sum_{j=1}^{k} \frac{\kappa}{Y_0} \Delta_n\right)^p \le \left(\frac{\kappa T}{Y_0}\right)^p.
\end{equation*}

Moreover, from Jensen's inequality,
\begin{equation*}
\left(\frac{\theta}{2}\sum_{j=1}^{k} \hat Y^n_{t^n_j} \Delta_n\right)^p \le \left(\frac{\theta}{2}\right)^p T^{p-1} \sum_{j=1}^{k} (\hat Y^n_{t^n_j})^p \Delta_n.
\end{equation*}

Finally,
\begin{equation*}
\left(\frac{\nu}{2}\right)^p |B^H_{t^n_k}|^p \le \left(\frac{\nu}{2}\right)^p \sup_{s\in[0,T]}|B^H_{s}|^p.
\end{equation*}

Hence, for all $t^n_k \le t^n_{\tau_1}$:
\begin{equation*}
\begin{gathered}
(\hat Y^n_{t^n_k})^p \le 4^{p-1} \left(Y_0^p + \left(\frac{\kappa T}{Y_0}\right)^p +  \left(\frac{\theta}{2}\right)^p T^{p-1} \sum_{j=1}^{k} (\hat Y^n_{t^n_j})^p \Delta_n + \left(\frac{\nu}{2}\right)^p \sup_{s\in[0,T]}|B^H_{s}|^p\right)
\\
\le \left((4Y_0)^p + \left(\frac{8 \kappa T}{Y_0}\right)^p + (8\nu)^p \sup_{s\in[0,T]}|B^H_s|^p\right) + (8\theta)^p T^{p-1}\sum_{j=1}^{k} (\hat Y^n_{t^n_j})^p \Delta_n.
\end{gathered}
\end{equation*}

\textbf{Step 2.} Assume that $\tau_1\ne n$, i.e. there are points of partition on the interval $(t^n_{\tau_1},T]$. From definition of $\tau_1$, $\hat Y^n_{t^n_{\tau_1}} \ge \frac{Y_0}{2}$ and for all points of the partition $t^n_k$ such that $t^n_k\in (t^n_{\tau_1},T]$:
\begin{equation*}
\left\{l=1,...,n~|~t^n_l \in (t^n_{\tau_1},t^n_k],~\hat Y^n_{t^n_l} < \frac{Y_0}{2}\right\}\ne \emptyset.
\end{equation*}

Let $t^n_k \in (t^n_{\tau_1},T]$ be fixed and denote
\begin{equation*}
\tau_2^k := \max\left\{l=1,...,n~|~t^n_l \in (t^n_{\tau_1},t^n_k],~\hat Y^n_{t^n_l} < \frac{Y_0}{2}\right\}.
\end{equation*}

It is obvious that $t^n_{\tau_1} < t^n_{\tau_2^k} \le t^n_k$ and $\hat Y^n_{t^n_{\tau_2^k}} < \frac{Y_0}{2}$, and
\begin{equation}\label{moments of approx: eq 1}
\begin{gathered}
(\hat Y^n_{t^n_k})^p = (\hat Y^n_{t^n_k}-\hat Y^n_{t^n_{\tau_2^k}}+\hat Y^n_{t^n_{\tau_2^k}})^p \le 2^{p-1} \left(|\hat Y^n_{t^n_k}-\hat Y^n_{t^n_{\tau_2^k}}|^p + (\hat Y^n_{t^n_{\tau_2^k}})^p\right)
\\
\le 2^{p-1} \left(|\hat Y^n_{t^n_k}-\hat Y^n_{t^n_{\tau_2^k}}|^p + \left(\frac{Y_0}{2}\right)^p\right) \le 2^{p-1} |\hat Y^n_{t^n_k}-\hat Y^n_{t^n_{\tau_2^k}}|^p + Y_0^p.
\end{gathered}
\end{equation}

In addition, if $t^n_{\tau_2^k} = t^n_k$, 
\begin{equation*}
|\hat Y^n_{t^n_k}-\hat Y^n_{t^n_{\tau_2^k}}|^p =  0,
\end{equation*}
and if $t^n_{\tau_2^k} < t^n_k$,
\begin{equation*}
\begin{gathered}
\left|\hat Y^n_{t^n_k}-\hat Y^n_{t^n_{\tau_2^k}}\right|^p = \left|\frac{1}{2}\sum_{j=\tau_2^k+1}^{k}\left(\frac{\kappa}{\hat Y^n_{t^n_j}} - \theta \hat Y^n_{t^n_j}\right) \Delta_n +\frac{\nu}{2}\left( B^H_{t^n_k} - B^H_{t^n_{\tau_2^k}} \right)\right|^p
\\
\le 4^{p-1} \left(\left(\frac{1}{2}\sum_{j=\tau_2^k+1}^{k} \frac{\kappa}{\hat Y^n_{t^n_j}} \Delta_n\right)^p + \left(\frac{\theta}{2}\sum_{j=\tau_2^k+1}^{k} \hat Y^n_{t^n_j} \Delta_n\right)^p + \left(\frac{\nu}{2}\right)^p |B^H_{t^n_k}|^p + \left(\frac{\nu}{2}\right)^p |B^H_{t^n_{\tau_2^k}}|^p\right).
\end{gathered}
\end{equation*}

From definition of $\tau_2^k$, for all points of the partition $t^n_l\in(t^n_{\tau_2^k},t^n_k]$ it holds that $\hat Y^n_{t^n_k}\ge \frac{Y_0}{2}$, so
\begin{equation*}
\left(\frac{1}{2}\sum_{j=\tau_2^k+1}^{k} \frac{\kappa}{\hat Y^n_{t^n_j}} \Delta_n\right)^p \le \left(\frac{\kappa T}{Y_0}\right)^p.
\end{equation*}

Furthermore,
\begin{equation*}
\left(\frac{\theta}{2}\sum_{j=\tau_2^k+1}^{k} \hat Y^n_{t^n_j} \Delta_n\right)^p \le \left(\frac{\theta}{2}\right)^p T^{p-1}\sum_{j=1}^{k} (\hat Y^n_{t^n_j})^p \Delta_n,
\end{equation*}
and
\begin{equation*}
\left(\frac{\nu}{2}\right)^p |B^H_{t^n_k}|^p + \left(\frac{\nu}{2}\right)^p |B^H_{t^n_{\tau_2^k}}|^p \le 2 \left(\frac{\nu}{2}\right)^p \sup_{s\in[0,T]}|B^H_s|^p.
\end{equation*}

Hence,
\begin{equation}\label{moments of approx: eq 2}
\begin{gathered}
\left|\hat Y^n_{t^n_k}-\hat Y^n_{t^n_{\tau_2^k}}\right|^p \le
\\
\le 4^{p-1} \left(\left(\frac{\kappa T}{Y_0}\right)^p + \left(\frac{\theta}{2}\right)^p T^{p-1}\sum_{j=1}^{k} (\hat Y^n_{t^n_j})^p \Delta_n + 2 \left(\frac{\nu}{2}\right)^p \sup_{s\in[0,T]}|B^H_s|^p\right).
\end{gathered}
\end{equation}

Finally, from \eqref{moments of approx: eq 1} and \eqref{moments of approx: eq 2},
\begin{equation*}
\begin{gathered}
(\hat Y^n_{t^n_k})^p \le 8^{p-1}\left(\left(\frac{\kappa T}{Y_0}\right)^p + \left(\frac{\theta}{2}\right)^p T^{p-1}\sum_{j=1}^{k} (\hat Y^n_{t^n_j})^p \Delta_n + 2 \left(\frac{\nu}{2}\right)^p \sup_{s\in[0,T]}|B^H_s|^p\right) + Y_0^p
\\
\le \left((4Y_0)^p + \left(\frac{8\kappa T}{Y_0}\right)^p + (8\nu)^p \sup_{s\in[0,T]}|B^H_s|^p\right) + (8\theta)^p T^{p-1}\sum_{j=1}^{k} (\hat Y^n_{t^n_j})^p \Delta_n.
\end{gathered}
\end{equation*}

Therefore, \eqref{pre-Gronwall} indeed holds for any point $t^n_k$ of the partition.

\textbf{Step 3.} As $n>2(8\theta)^pT^{p-1}$, 
\[
\frac{1}{2} \le 1-(8\theta)^p T^{p-1}\Delta_n \le 1,
\]
therefore, as, due to \eqref{pre-Gronwall},
\begin{equation*}
\begin{gathered}
\left(1-(8\theta)^p T^{p-1}\Delta_n\right) (\hat Y^n_{t^n_k})^p 
\\
\le \left((4Y_0)^p + \left(\frac{8\kappa T}{Y_0}\right)^p + (8\nu)^p \sup_{s\in[0,T]}|B^H_s|^p\right) + (8\theta)^p T^{p-1}\sum_{j=1}^{k-1} (\hat Y^n_{t^n_j})^p \Delta_n,
\end{gathered}
\end{equation*}
we have
\begin{equation*}
\begin{gathered}
(\hat Y^n_{t^n_k})^p \le 2\left((4Y_0)^p + \left(\frac{8\kappa T}{Y_0}\right)^p + (8\nu)^p \sup_{s\in[0,T]}|B^H_s|^p\right) + 2(8\theta)^p T^{p-1}\sum_{j=1}^{k-1} (\hat Y^n_{t^n_j})^p \Delta_n.
\end{gathered}
\end{equation*}

Using the discrete version of the Gr\"onwall's lemma, we obtain:
\begin{equation*}
(\hat Y^n_{t^n_k})^p \le 2\left((4Y_0)^p + \left(\frac{8\kappa T}{Y_0}\right)^p + (8\nu)^p \sup_{s\in[0,T]}|B^H_s|^p\right)e^{2(8\theta T)^p},
\end{equation*}
i.e., taking into account that the right-hand side does not depend on $n$ and remarks in the beginning of the proof, there is such $C>0$ that
\begin{equation}\label{upper boundary}
\sup_{n\ge 0}\sup_{t\in[0,T]}(\hat Y^n_{t})^p < C(1+\sup_{t\in[0,T]}|B^H_t|^p).
\end{equation}

Now the claim of the Theorem follows from the fact that the right-hand side of \eqref{upper boundary} does not depend on $n$ and that (see, for example, \cite{MishuraFBM})
\[
\E \sup_{s\in[0,T]}|B^H_s|^p < \infty.
\]
\QED
\\[11pt]
\textbf{Proof of Corollary \ref{exp moments of approximations}.} From \eqref{upper boundary} it follows that there is such $C>0$ that
\[
\sup_{n\ge 0}\sup_{t\in[0,T]}(\hat Y^n_{t})^\varrho < C(1+\sup_{t\in[0,T]}|B^H_t|^\varrho).
\]
The rest of the proof is similar to Theorem \ref{solution to price equation}, 1. \QED
\\[11pt]
\textbf{Proof of Theorem \ref{approx conv theor}.} We shall proceed as in proof of Lemma 14, \cite{BdPM}.

Using H\"older's inequality, we write:
\begin{equation*}
\begin{gathered}
\E|X_T - \hat X^n_T|^2 = \E\left| - \frac{1}{2}\int_0^T \sigma^2(Y_s)ds+\int_0^T \sigma(Y_s)dW_s + \frac{1}{2} \int_0^T \sigma^2(\hat Y^n_{s})ds -  \int_0^T \sigma(\hat Y^n_{s})dW_s\right|^2
\\
\le C\left( \E\left| - \frac{1}{2}\int_0^T (\sigma^2(Y_s) - \sigma^2(\hat Y^n_{s})) ds\right|^2 + \E\left| \int_0^T (\sigma( Y_{s}) -\sigma(\hat Y^n_{s}) )dW_s\right|^2\right)
\\
\le C\left( \int_0^T \E [\sigma^2(Y_s) - \sigma^2(\hat Y^n_{s})]^2 ds + \int_0^T \E[\sigma( Y_{s}) -\sigma(\hat Y^n_{s}) ]^2ds\right)
\\
= C\left( \int_0^T \E [(\sigma(Y_s) - \sigma(\hat Y^n_{s}))(\sigma(Y_s) + \sigma(\hat Y^n_{s}))]^2 ds + \int_0^T \E[\sigma( Y_{s}) -\sigma(\hat Y^n_{s}) ]^2ds\right).
\end{gathered}
\end{equation*}

From Assumption \ref{AssumSigma} (iii), Jensen's inequality and Theorem \ref{China theorem},
\begin{equation*}
\begin{aligned}
\int_0^T \E[\sigma( Y_{s}) -\sigma(\hat Y^n_{s}) ]^2ds &\le C_\sigma^2 \int_0^T \E[(Y_{s}-\hat Y^n_{s})^{2r} ]ds 
\\
&\le C_\sigma^2 \int_0^T \left(\E[(Y_{s}-\hat Y^n_{s})^{4} ]\right)^{\frac{r}{2}} ds
\\
& \le C \Delta_n^{2rH}.
\end{aligned}
\end{equation*}

Moreover, Assumption \ref{AssumSigma}, (ii) and (iii), implies that
\begin{equation*}
\begin{gathered}
\E \left[\left(\sigma(Y_s) - \sigma(\hat Y^n_{s})\right)\left(\sigma(Y_s) + \sigma(\hat Y^n_{s})\right)\right]^2 \le C^2_\sigma\E \left[(Y_s - \hat Y^n_{s})^{2r}\left(2\sigma^2(Y_s) + 2\sigma^2(\hat Y^n_{s})\right)\right]
\\
\le C \E \left[(Y_s - \hat Y^n_{s})^{2r}\left((1+Y^q_s)^2 + (1+(\hat Y^n_{s})^q)^2\right)\right]
\\
\le C\E \left[(Y_s - \hat Y^n_{s})^{2r}\left(1+Y^{2q}_s + (\hat Y^n_{s})^{2q}\right)\right]
\\
\le C \left(\E(Y_s-\hat Y^n_s)^{4r}\right)^{\frac{1}{2}} \left(\E\left[1 + Y^{4q}_s + (\hat Y^n_{s})^{4q}\right]\right)^{\frac{1}{2}}.
\end{gathered}
\end{equation*}

From Theorem \ref{China theorem},
\begin{equation*}
\left(\E(Y_s-\hat Y^n_s)^{4r}\right)^{\frac{1}{2}} \le \left(\E(Y_s-\hat Y^n_s)^{4}\right)^{\frac{r}{2}} \le C\Delta_n^{2rH},
\end{equation*}
and, from Theorems \ref{UpperBoundForY} and \ref{th moments of approximations},
\begin{equation*}
\left(\E\left[1 + Y^{4q}_s + (\hat Y^n_{s})^{4q}\right]\right)^{1/2} < \infty.
\end{equation*}

Therefore, taking into account bounds above, there is such constant $C>0$ that
\[
\E|X_T - \hat X^n_T|^2 \le C \Delta^{2rH}.
\]

Now, let us prove \eqref{Zapprox}. Taking into account Assumption \ref{AssumSigma} (i), 
\[
\left|\frac{1}{\sigma(x)} - \frac{1}{\sigma(y)}\right| = \frac{|\sigma(x) - \sigma(y)|}{\sigma(x)\sigma(y)} \le \frac{|\sigma(x) - \sigma(y)|}{\sigma^2_{\min}},
\]
so, from Assumption \ref{AssumSigma} (iii),
\begin{equation*}
\begin{aligned}
\E(Z_T - \hat Z^n_T)^2 &= \int_0^T\E\left(\frac{1}{\sigma(Y_s)} - \frac{1}{\sigma(\hat Y^n_s)}\right)^2ds 
\\
&\le \frac{1}{\sigma^2_{\min}} C_\sigma \int_0^T \E(Y_s - \hat Y^n_s)^{2r}ds
\\
& \le C \int_0^T \left(\E(Y_s - \hat Y^n_s)^{4}\right)^{\frac{r}{2}}ds
\\
& \le C\Delta_n^{2rH}.
\end{aligned}
\end{equation*}
\QED
\\[11pt]
\textbf{Proof of Lemma \ref{auxiliary approx lemma}.} It is clear that 
\begin{equation}\label{I1 and I2}
\E\left|\frac{F(S_T)}{S_T} - \frac{F(\hat S^n_T)}{\hat S^n_T}\right|^2 \le 2\E\left|\frac{F(S_T)}{S_T} - \frac{F(S_T)}{\hat S^n_T}\right|^2 + 2\E\left|\frac{F(S_T)}{\hat S^n_T} - \frac{F(\hat S^n_T)}{\hat S^n_T}\right|^2.
\end{equation}

Now we shall estimate the right-hand side of \eqref{I1 and I2} term by term.

\begin{equation*}
\E \left|F(S_T)\left(\frac{1}{S_T} - \frac{1}{\hat S^n_T}\right)\right|^2 \le \left(\E\left(F(S_T)\right)^4 \E\left(\frac{1}{S_T} - \frac{1}{\hat S^n_T}\right)^4\right)^{\frac{1}{2}}.
\end{equation*}

From Assumption \ref{AssumPayoff} (i), both $f$ and $F$ are of polynomial growth, therefore, due to \eqref{UpperBoundForSn},
\[
\E\left(F(S_T)\right)^4 < \infty.
\]

Furthermore, using sequentially the inequalities
\begin{equation*}
\begin{gathered}
|e^x - e^y| \le (e^x + e^y)|x-y|, \quad x, y \in \mathbb R,
\\
(x+y)^{2n} \le C(n)(x^{2n} + y^{2n}), \quad x, y \in \mathbb R, \quad n\in\mathbb N.
\end{gathered}
\end{equation*}
and H\"older's inequality, we obtain that
\begin{equation*}
\begin{gathered}
\E\left(\frac{1}{S_T} - \frac{1}{\hat S^n_T}\right)^4 = \frac{1}{S_0^4 e^{4\mu t}}\E\left[e^{\frac{1}{2}\int_0^T\sigma^2(Y_s)ds - \int_0^T \sigma (Y_s)dW_s} - e^{\frac{1}{2}\int_0^T\sigma^2(\hat Y^n_s)ds - \int_0^T \sigma (\hat Y^n_s)dW_s}\right]^4
\\
\le C \E\Bigg[\left(e^{\frac{1}{2}\int_0^T\sigma^2(Y_s)ds - \int_0^T \sigma (Y_s)dW_s} + e^{\frac{1}{2}\int_0^T\sigma^2(\hat Y^n_s)ds - \int_0^T \sigma (\hat Y^n_s)dW_s}\right)^4
\\
\times\bigg(\frac{1}{2}\int_0^T\sigma^2(Y_s)ds - \int_0^T \sigma (Y_s)dW_s - \frac{1}{2}\int_0^T\sigma^2(\hat Y^n_s)ds + \int_0^T \sigma (\hat Y^n_s)dW_s\bigg)^4\Bigg]
\\
\le C \E\Bigg[\left(e^{2\int_0^T\sigma^2(Y_s)ds - 4\int_0^T \sigma (Y_s)dW_s} + e^{2\int_0^T\sigma^2(\hat Y^n_s)ds - 4\int_0^T \sigma (\hat Y^n_s)dW_s}\right)
\\
\times\bigg(\frac{1}{2}\int_0^T\sigma^2(Y_s)ds - \int_0^T \sigma (Y_s)dW_s - \frac{1}{2}\int_0^T\sigma^2(\hat Y^n_s)ds + \int_0^T \sigma (\hat Y^n_s)dW_s\bigg)^4\Bigg]
\\
\le C \left(\E\left[e^{4\int_0^T\sigma^2(Y_s)ds - 8\int_0^T \sigma (Y_s)dW_s} + e^{4\int_0^T\sigma^2(\hat Y^n_s)ds - 8\int_0^T \sigma (\hat Y^n_s)dW_s}\right]\right)^{\frac{1}{2}}
\\
\times \left(\E\left[\bigg(\frac{1}{2}\int_0^T\sigma^2(Y_s)ds - \int_0^T \sigma (Y_s)dW_s - \frac{1}{2}\int_0^T\sigma^2(\hat Y^n_s)ds + \int_0^T \sigma (\hat Y^n_s)dW_s\bigg)^8\right]\right)^{\frac{1}{2}}.
\end{gathered}
\end{equation*}

Next, from \eqref{DDexpY} and Remark \ref{moments of approximations}  it follows that
\begin{equation*}
\begin{gathered}
\E\left[e^{4\int_0^T\sigma^2(Y_s)ds - 8\int_0^T \sigma (Y_s)dW_s}\right] < \infty,
\\
\E\left[ e^{4\int_0^T\sigma^2(\hat Y^n_s)ds - 8\int_0^T \sigma (\hat Y^n_s)dW_s}\right] < \infty,
\end{gathered}
\end{equation*}
so, using this together with H\"older and Burkholder-Davis-Gundy inequalities, we continue the chain as follows:
\begin{equation*}
\begin{gathered}
\E\left(\frac{1}{S_T} - \frac{1}{\hat S^n_T}\right)^4
\\
\le C \left(\E\left[\bigg(\frac{1}{2}\int_0^T\sigma^2(Y_s)ds - \int_0^T \sigma (Y_s)dW_s - \frac{1}{2}\int_0^T\sigma^2(\hat Y^n_s)ds + \int_0^T \sigma (\hat Y^n_s)dW_s\bigg)^8\right]\right)^{\frac{1}{2}}
\\
\le C \left(\E\left[\bigg(\int_0^T\sigma^2(Y_s)ds - \int_0^T\sigma^2(\hat Y^n_s)ds\bigg)^8\right] + \E\left[\bigg(\int_0^T \sigma (Y_s)dW_s - \int_0^T \sigma (\hat Y^n_s)dW_s\bigg)^8\right]\right)^{\frac{1}{2}}
\\
\le C\left(\E\left[\int_0^T\left(\sigma^2(Y_s) - \sigma^2(\hat Y^n_s)\right)^8ds\right] + \E\left[\left(\int_0^T\left(\sigma(Y_s) - \sigma(\hat Y^n_s)\right)^2 ds\right)^4\right]\right)^{\frac{1}{2}}
\\
\le C \Bigg(\int_0^T\E\left(\left(\sigma(Y_s) - \sigma(\hat Y^n_s)\right)\left(\sigma(Y_s) + \sigma(\hat Y^n_s)\right)\right)^8ds
\\
+ \int_0^T\E\left(\sigma(Y_s) - \sigma(\hat Y^n_s)\right)^8 ds\Bigg)^{\frac{1}{2}}.
\end{gathered}
\end{equation*}

By applying Assumption \ref{AssumSigma}, (ii) and (iii),
\begin{equation*}
\begin{gathered}
\int_0^T\E\left(\left(\sigma(Y_s) - \sigma(\hat Y^n_s)\right)\left(\sigma(Y_s) + \sigma(\hat Y^n_s)\right)\right)^8ds
\\
\le C\int_0^T \E\left[\left(Y_s - \hat Y^n_s\right)^{8r}\left(1+ Y_s^q + (\hat Y^n_s)^q\right)^8\right]ds
\\
\le C \int_0^T \left(\E\left[\left(Y_s - \hat Y^n_s\right)^{16r}\right]\right)^{\frac{1}{2}}\left(\E\left[1+ Y_s^{16q} + (\hat Y^n_s)^{16q}\right]\right)^{\frac{1}{2}}ds
\end{gathered}
\end{equation*}
and
\begin{equation*}
\begin{gathered}
 \int_0^T\E\left(\sigma(Y_s) - \sigma(\hat Y^n_s)\right)^8ds \le C\int_0^T\E\left(Y_s - \hat Y^n_s\right)^{8r}ds.
\end{gathered}
\end{equation*}

From Theorems \ref{UpperBoundForY} and \ref{th moments of approximations},
\begin{equation*}
\begin{gathered}
\sup_{n\ge 1}\sup_{t\in[0,T]} \E\left[1+ Y_s^{16q} + (\hat Y^n_s)^{16q}\right] < \infty,
\end{gathered}
\end{equation*}
and, according from Theorem \ref{China theorem},
\begin{equation*}
\begin{gathered}
\int_0^T\left(\E\left[\left(Y_s - \hat Y^n_s\right)^{16r}\right]\right)^{\frac{1}{2}}ds \le C\Delta_n^{16rH}
\\
\int_0^T\E\left(Y_s - \hat Y^n_s\right)^{8r}ds \le C \Delta_n^{8rH},
\end{gathered}
\end{equation*}
hence
\begin{equation}\label{I1}
\begin{gathered}
\E\left(\frac{1}{S_T} - \frac{1}{\hat S^n_T}\right)^4 \le C\Delta_n^{4rH} \le C \Delta_n^{2rH}.
\end{gathered}
\end{equation}

Now, let us move to the second term of the right-hand side of \eqref{I1 and I2}.
\begin{equation*}
\begin{gathered}
\E\left|\frac{F(S_T)}{\hat S^n_T} - \frac{F(\hat S^n_T)}{\hat S^n_T}\right|^2 \le \left(\E\left[\frac{1}{(\hat S^n_T)^4}\right]\right)^{\frac{1}{2}}\left(\E\left(F(S_T)-F(\hat S^n_T)\right)^{4}\right)^{\frac{1}{2}}.
\end{gathered}
\end{equation*}

Due to Remark \ref{moments of approximations}, 
\[
\E\left[\frac{1}{(\hat S^n_T)^4}\right] < \infty,
\]
and, from Assumption \ref{AssumPayoff} (i),
\begin{equation*}
\begin{gathered}
\left(\E\left(F(S_T)-F(\hat S^n_T)\right)^{4}\right)^{\frac{1}{2}} = \left(\E\left(\int_{S_T\wedge \hat S^n_T}^{S_T\vee \hat S^n_T}f(x)dx\right)^4\right)^{\frac{1}{2}}
\\
\le C\left(\E\left[(S_T - \hat S^n_T)^4 (1 + S_T^p + (\hat S^n_T)^p)^4\right]\right)^{\frac{1}{2}} \le C \left(\E((S_T - \hat S^n_T)^8 \E(1 + S_T^p + (\hat S^n_T)^p)^8\right)^{\frac{1}{4}}.
\end{gathered}
\end{equation*}

According to \eqref{UpperBoundForSn} and Remark \ref{moments of approximations},
\[
\E(1 + S_T^p + (\hat S^n_T)^p)^8 < \infty,
\]
so
\[
\E\left|\frac{F(S_T)}{\hat S^n_T} - \frac{F(\hat S^n_T)}{\hat S^n_T}\right|^2 \le C \left(\E((S_T - \hat S^n_T)^8\right)^{\frac{1}{4}}.
\]

To get the final result, we can proceed just as in the upper bound for the first term in the right-hand side of \eqref{I1 and I2}. Thus
\begin{equation}\label{I2}
\begin{gathered}
\E\left|\frac{F(S_T)}{\hat S^n_T} - \frac{F(\hat S^n_T)}{\hat S^n_T}\right|^2 \le C \left(\int_0^T \left(\left(\E|Y_s - \hat Y^n_s|^{32r}\right)^{\frac{1}{2}} + \E|Y_s - \hat Y^n_s|^{16r}\right)ds\right)^{\frac{1}{8}}
\\
\le C\Delta_n^{2rH}.
\end{gathered}
\end{equation}

Relations \eqref{I1} and \eqref{I2} together with \eqref{I1 and I2} complete the proof. \QED
\\[11pt]
\textbf{Proof of Theorem \ref{final option pricing}.} According to Theorem \ref{OptionPriceFormula},
\begin{equation*}
\begin{gathered}
\left|\E f(S_T) - \E\left[\frac{F(\hat S^n_T)}{\hat S^n_T}\left(1 + \frac{\hat Z^n_T}{T}\right)\right]\right| = \left|\E\left[\frac{F(S_T)}{S_T}\left(1 + \frac{Z_T}{T}\right)\right] - \E\left[\frac{F(\hat S^n_T)}{\hat S^n_T}\left(1 + \frac{\hat Z^n_T}{T}\right)\right]\right|
\\
\le \frac{1}{T}\E\left[\left|\frac{F(S_T)}{S_T}(Z_T - \hat Z^n_T)\right|\right] + \E\left[\left|\left(1+\frac{\hat Z^n_T}{T}\right)\left(\frac{F(S_T)}{S_T} - \frac{F(\hat S^n_T)}{\hat S^n_T}\right)\right|\right]
\\
\le \frac{1}{T} \left(\E\left(\frac{F(S_T)}{S_T}\right)^2 \E(Z_T - \hat Z_T^n)^2\right)^{\frac{1}{2}} + \left(\E\left(1+\frac{\hat Z^n_T}{T}\right)^2 \E\left(\frac{F(S_T)}{S_T} - \frac{F(\hat S^n_T)}{\hat S^n_T}\right)^2\right)^{\frac{1}{2}}.
\end{gathered}
\end{equation*}

According to Theorem \ref{solution to price equation}, Assumption \ref{AssumPayoff} (i) and the Cauchy-Schwartz inequality, $\E\left(\frac{F(S_T)}{S_T}\right)^2<\infty$. Next,
\[
\sup_{n\ge 1} \E(Z_T^n)^2 < \frac{T}{\sigma^2_{\min}}.
\]

The proof now follows from Theorem \ref{approx conv theor} and Lemma \ref{auxiliary approx lemma}. \QED

\appendix
\section{Necessary results from Malliavin Calculus}\label{sec: preliminaries}

In this section, we recall several main definitions and results related to Malliavin calculus. For more detail, we refer to \cite{Nual}.

Let $B^H = \{B^H_t,~t\in[0,T]\}$ be a fractional Brownian motion with $H \in [1/2,1)$ on the standard probability space $\{\Omega, \mathcal F, \mathbb F = \{\mathcal F_t\}_{t\in[0,T]}, \mathbb P\}$, where $\Omega = C([0,T],\mathbb R)$, i.e. a centered Gaussian process that starts in zero and has a covariance function of the form
\begin{equation*}
R_H(t,s) := \E B^H_t B^H_s = \frac{1}{2} (t^{2H} + s^{2H} -|t-s|^{2H}), \quad s,t\in[0,T].
\end{equation*}

Note that the covariance function of the fractional Brownian motion has the form 
\begin{equation*}
R_H(t,s) = \begin{cases} \int_0^T \mathbbm 1_{[0,t]}(u) \mathbbm 1_{[0,s]}(u)du, \quad H=\frac{1}{2},
\\
\int_0^t \int_0^s \varphi (\tau, u)du d\tau, \quad H>\frac{1}{2},
\end{cases}
\end{equation*}
where $\varphi(\tau,u) := H(2H-1)|u-\tau|^{2H-2}$. 

On the set of all step functions on $[0,T]$, define an inner product that acts as follows for the indicator functions:
\begin{equation*}
\langle \mathbbm 1_{[0,t]}, \mathbbm 1_{[0,s]} \rangle_{\mathcal H} := R_H(t,s). 
\end{equation*}

Denote $\mathcal H$ the Hilbert space that is the closure of the space of all step functions on $[0,T]$ with respect to $\langle \cdot, \cdot \rangle_{\mathcal H}$.

\begin{remark}\label{HilbWiener}
If $H=1/2$, $\mathcal H$ coincides with $L_2([0,T])$.
\end{remark}

The mapping $\mathbbm 1_{[0,t]} \to B^H_t$ can be extended to a linear isometry from $\mathcal H$ onto a closed subspace $\mathcal H_1$ of $L^2(\Omega,\mathcal F,\mathcal P)$ associated with $B^H$. We will denote this isometry by $\phi \to B^H_\phi$. In this case, for all $\phi, \psi \in \mathcal H$:
\begin{equation*}
\langle \phi, \psi \rangle_{\mathcal H} = \E B^H_\phi B^H_\psi.
\end{equation*}

Denote by $C^{\infty}_p(\mathbb R^n)$ the set of all infinitely differentiable functions with the derivatives of at most polynomial growth at infinity.

\begin{definition}
Random variables $\xi$ of the form 
\begin{equation}\label{smooth rv}
\xi = h(B^H_{\phi_1}, ..., B^H_{\phi_n}),
\end{equation}
where $h\in C^{\infty}_p(\mathbb R^n)$, $\phi_1,...,\phi_n\in \mathcal H$, $n\ge 1$, are called smooth.
\end{definition}

Denote $\mathcal S$ the set of all smooth random variables.

\begin{definition}
Let $\xi \in \mathcal S$. The stochastic or Malliavin derivative of a smooth random variable $\xi$ of the form \eqref{smooth rv} is the $\mathcal H$-valued random variable given by
\begin{equation*}
D\xi = \sum_{i=1}^n \frac{\partial h}{\partial x_i} (B^H_{\phi_1}, ..., B^H_{\phi_n}) \phi_i.
\end{equation*}
\end{definition}

\begin{remark}
If $\phi_i = \mathbbm 1_{[0,t_i]}$, $t_i\in[0,T]$, $i=1,...,n$, then $B^H_{\mathbbm 1_{[0,t_i]}} = B^H_{t_i}$ and the real-valued random variable of the form
\begin{equation*}
D_t \xi = \sum_{i=1}^n \frac{\partial h}{\partial x_i} (B^H_{t_1}, ..., B^H_{t_n}) \mathbbm 1_{[0,t_i]}(t), \quad t\in[0,T],
\end{equation*} 
is called the stochastic derivative of $\xi$ at time $t$.
\end{remark}

According to Proposition 1.2.1 from \cite{Nual}, $D$ as an operator from the subset of $L^p(\Omega)$ to $L^p(\Omega,\mathcal H)$ is closable for any $p\ge 1$ and we shall use the same notation $D$ for the closure. 
\begin{definition}
Let $p\ge1$. The domain $\mathbb D^{1,p}$ of $D$ is the closure of the class of smooth random variables $\mathcal S$ with respect to the norm
\begin{equation*}
\lVert \xi \rVert_{1,p} := \left(\E|\xi|^p + \E\lVert D\xi\rVert^p_{\mathcal H}\right)^{1/p}.
\end{equation*}
\end{definition}

\begin{remark}
For $p=2$, the space $\mathbb D^{1,2}$ is the Hilbert space with respect to the inner product
\begin{equation*}
\langle \xi, \eta \rangle_{1,2} = \E \xi\eta + \E\left[\langle D\xi, D\eta\rangle_{\mathcal H}\right].
\end{equation*}
\end{remark}

\begin{proposition}{(\cite{Nual}, Proposition 1.2.3)} 
Let $F$: $\mathbb R^m \to \mathbb R$ be a continuously differentiable function with bounded partial derivatives, and fix $p\ge 1$. Suppose that $\xi = (\xi_1,...,\xi_m)$ is a random vector whose components belong to the space $\mathbb D^{1,p}$. Then $F(\xi) \in \mathbb D^{1,p}$ and
\begin{equation*}
DF(\xi) = \sum_{i=1}^m \frac{\partial F(\xi)}{\partial x_i}  D\xi_i.
\end{equation*}
\end{proposition}

\begin{remark}
In what follows, we will consider the case $p=2$.
\end{remark}

\begin{definition}
The divergence or Skorokhod operator $\delta$ is the adjoint of the operator $D$, i.e. an undounded operator on $L^2(\Omega,\mathcal H)$ with values in $L^2(\Omega)$ such that:
\begin{itemize}
\item[(i)] the domain of $\delta$, denoted by $\mathsf{Dom~}\delta$, is the set of $\mathcal H$-valued square integrable random variables $\zeta\in L^2(\Omega,\mathcal H)$ such that for all $\xi\in\mathbb D^{1,2}$:
\begin{equation*}
\left|\E\left[\langle D\xi, \zeta\rangle_{\mathcal H}\right]\right| \le C_\zeta (\E\xi^2)^{1/2},
\end{equation*}
where $C_\zeta$ is some constant depending on $\zeta$;
\item[(ii)] if $\zeta$ belongs to  $\mathsf{Dom~}\delta$, then $\delta(\zeta)$ is the element of $L^2(\Omega)$ characterized by
\begin{equation*}
\E [\xi \delta(\zeta)] = \E \left[\langle D\xi,\zeta\rangle_{\mathcal H}\right]
\end{equation*}
for any $\xi\in\mathbb D^{1,2}$.
\end{itemize}
\end{definition}

The Skorokhod operator $\delta$ is closed.

\begin{remark} Let $V = \{V_t,~t\in[0,T]\}$ be the Wiener process, $\mathcal H_V = L_2([0,T])$ be the associated Hilbert space (see Remark \ref{HilbWiener}) and $\delta_V$ be the corresponding divergence operator. In this case, the elements of $\mathsf{Dom}~\delta_V \subset L^2([0,T]\times\Omega)$ are square-integrable processes, and the divergence $\delta_V(\zeta)$ is called the Skorokhod stochastic integral of the process $\zeta$ with respect to $V$ and is denoted as follows:
\begin{equation*}
\delta_V(\zeta) = \int_0^T \zeta_t \delta V_t.
\end{equation*}
 According to \cite{Nual}, Section 1.3.2, the Skorokhod integral is correctly defined for all elements of the space $\mathbb L^{1,2} = L^2([0,T], \mathbb D^{1,2})$ with the norm $\lVert\cdot\rVert_{\mathbb L^{1,2}}$ such that
\begin{equation*}
\lVert \zeta \rVert_{\mathbb L^{1,2}} = \E\left(\int_0^T \zeta^2_tdt + \int_0^T\int_0^T (D_s\zeta_t)^2dtds\right).
\end{equation*}
\end{remark}

\begin{remark}
Let $B^H$ be a fractional Brownian motion with $H>1/2$. Similarly to the Wiener process case, we shall call the corresponding divergence $\delta_H(\zeta)$ the Skorokhod stochastic integral with respect to fractional Brownian motion and shall denote it as
\[
\delta_H(\zeta) = \int_0^T \zeta_t \delta B^H_t.
\]
\end{remark}

In what follows, we shall use the definition of pathwise stochastic integral with respect to fractional Brownian motion proposed in \cite{Zahle} and denote it by $\int_0^T \zeta_t dB^H_t$. There is a useful result that connects stochastic and Skorokhod integrals, which is given below.

Let $H>1/2$ and 
\[
|\mathcal H| = \left\{\phi\in\mathcal H~\bigg|~\lVert \phi\rVert_{|\mathcal H|}^2 = \int_0^T\int_0^T |\phi(\tau)||\phi(u)|\varphi(\tau, u)dud\tau < \infty\right\}.
\]

\begin{theorem}[\cite{Nual}, Proposition 5.2.1]\label{NualTheorem} Let $\zeta = \{\zeta_t,~t\in[0,T]\}$ be a stochastic process in the space $\mathbb D^{1,2}(|\mathcal H|)$ with H\"older continuous trajectories up to the order $H$ and $D_s^H$ be the Malliavin derivative operator with respect to $B^H$. Suppose that a.s.
\[
\int_0^T\int_0^T |D^H_s u_t| |t-s|^{2H-2} dsdt < \infty.
\]
Then $\zeta$ is Stratonovich integrable and
\[
\int_0^T \zeta_t \circ dB^H_t = \int_0^T u_t \delta B^H_t + \int_0^T \int_0^T D^H_s \zeta_t \varphi(s,t) dsdt.
\]
\end{theorem}

\section*{Acknowledgments}

The first author acknowledges that the present research is carried through within the frame and support of the ToppForsk project nr. 274410 of the Research Council of Norway with title STORM: Stochastics for Time-Space Risk Models. The second author was partially supported by the grant 346300 for IMPAN from the Simons Foundation and the matching 2015-2019 Polish MNiSW fund.

\end{document}